\documentclass[12pt]{amsart}
\usepackage{amsmath,amssymb,amscd,array, mathrsfs }

\usepackage{amsmath,amscd,amssymb,amsthm,array}
\usepackage{color}
\usepackage{cancel}
\usepackage{amsmath,amssymb,mathrsfs,amsthm, tikz-cd,mathrsfs}
\usepackage{comment}
\def\url#1{\expandafter\s

\tring\csname #1\endcsname}

\def\mmat #1,#2,#3,#4,{\text{\small\arraycolsep=3pt $
\begin{pmatrix}#1&#2\\#3&#4\end{pmatrix}$}}

\usepackage{hyperref}

\usepackage{comments}
\newComments\SBe{Said}{blue}
\newComments\SBo{Sofiane}{blue}
\newComments\AM{Nacer}{blue}
\newComments\DL{DL}{red}
\newComments\QEh{QEh}{blue}

\def\mmat #1,#2,#3,#4,{\text{\small\arraycolsep=3pt $
\begin{pmatrix}#1&#2\\#3&#4\end{pmatrix}$}}

\usepackage{lscape}
\usepackage{tikz-cd}

\usepackage{DLdef1}
\usepackage{multicol}

\def\mmat #1,#2,#3,#4,{\text{\small\arraycolsep=3pt $
\begin{pmatrix}#1&#2\\#3&#4\end{pmatrix}$}}

\renewcommand {\ssbegin}[2][*]
 {\refstepcounter{subsection}%
\if#1*
\addcontentsline{toc}{subsection}{\thesubsection.\hskip 1pc #2}%
\else
\addcontentsline{toc}{subsection}{\thesubsection.\hskip 1pc #2. #1}%
\fi
 \def \secno {\gdef \secno {}{\ssecfont
\thesubsection.\hskip 2ex}%
 }%
 \begin{#2}}

\renewcommand {\sssbegin}[2][*]
 {\refstepcounter{subsubsection}
\if#1*
\addcontentsline{toc}{subsubsection}{\thesubsubsection.\hskip 1pc #2}%
\else
\addcontentsline{toc}{subsubsection}{\thesubsubsection.\hskip 1pc #2. #1}
\fi
 \def \secno {\gdef \secno {}{\ssecfont \thesubsubsection.\hskip 2ex}%
 }%
 \begin{#2}}

\renewcommand {\parbegin}[2][*]
 {\refstepcounter{paragraph}
\if#1*
\addcontentsline{toc}{paragraph}{\theparagraph.\hskip 1pc #2}%
\else
\addcontentsline{toc}{paragraph}{\theparagraph.\hskip 1pc #2. #1}
\fi
 \def \secno {\gdef \secno {}{\ssecfont \theparagraph.\hskip 2ex}%
 }%
 \begin{#2}}

\newcommand {\ce}{{\text{CE}}}
\newcommand {\obs}{{\text{obs}}}
\newcommand{\trr}{\triangleright}

\newcommand{\del}{\partial}
\rmnameii{etr}{etr}
\rmnameii{evv}{ev}

\DeclareMathOperator{\K}{\mathbb{K}}

\newcommand{\Z}{\mathbb{Z}}
\newcommand{\D}{\mathrm{d}}
\newcommand{\w}{\omega}

\begin{document}

\title[On Poisson superalgebras in characteristic 2]{On Poisson superalgebras in characteristic 2}



\author{Quentin Ehret}
\address {Division of Science and Mathematics, New York University Abu Dhabi, P.O. Box 129188, Abu Dhabi, United Arab Emirates.}
\email{qe209@nyu.edu}

\thanks{The author was supported by the grant NYUAD-065.}



\begin{abstract}
This paper is devoted to the study of Poisson superalgebras over fields of characteristic $2$. We investigate their representations, semidirect products, cohomology, formal deformations, and universal enveloping algebras. We also introduce Lie-Rinehart superalgebras in characteristic $2$ and clarify their connections with Poisson superalgebras. In particular, we show that the universal enveloping algebra of a Poisson superalgebra coincides with the universal enveloping algebra of an associated Lie-Rinehart superalgebra. We also compute examples and initiate the study of pre-Poisson superalgebras.

\end{abstract}


\maketitle

\thispagestyle{empty}
\setcounter{tocdepth}{2}
\tableofcontents

\section{Introduction} \label{sec:intro}
\subsection{Poisson (super)algebras in characteristic zero} Poisson brackets occur naturally in many different contexts and have numerous applications, particularly in classical mechanics. Poisson algebras were introduced as an algebraic framework for Poisson brackets and have since played an essential role in areas such as Hamiltonian mechanics, quantum mechanics, integrable systems, quantization, symplectic and differential geometry, see \cite{LPV} for a thorough survey of these structures. Roughly speaking, a Poisson algebra is an associative commutative algebra equipped with a Lie bracket that is a biderivation for the product. Poisson (co)homology was introduced by Lichnerowicz (\cite{Li}) and Koszul (\cite{K}) in the setting of Poisson manifolds. In the purely algebraic context, the construction of the Poisson cohomology relies heavily on the Chevalley–Eilenberg cohomology of Lie algebras. The usual Chevalley–Eilenberg cochain space is replaced by the space of skew-symmetric bilinear biderivations of the associative commutative product (see, for instance, \cite[Section~4.1.2]{LPV}). From the point of view of representation theory, Oh defined the universal enveloping algebra of a Poisson algebra in \cite{Oh}; see also \cite{OPS, LOV}. Poisson algebras are related to Lie-Rinehart algebras; this connection is explored in \cite{Hu}. A Lie-Rinehart algebra is a triple $(A, L, \rho)$, where $A$ is an associative commutative algebra, $L$ is both a Lie algebra and an $A$-module, and $\rho : L \to \Der(A)$ is a morphism of Lie algebras and $A$-modules, called \emph{anchor map}. This map satisfies a compatibility condition that controls, in particular, the failure of the Lie bracket on $L$ to be $A$-linear. Lie-Rinehart algebras first appeared in the works of Herz (\cite{He}) and Palais (\cite{P}), and were developed systematically by Rinehart (\cite{Ri}), in order to provide a formalism of differential forms for general commutative algebras. He also studied representations and cohomology of Lie-Rinehart algebras, which were later thoroughly studied by Huebschmann (see \cite{Hu}). In particular, Huebschmann showed that any Poisson algebra gives rise to a Lie-Rinehart algebra by considering its module of K\"{a}hler differentials, and he used this construction to obtain Poisson (co)homology in terms of the associated Lie-Rinehart algebra. Furthermore, the same construction provides a way to define the universal enveloping algebra of a Poisson algebra via the enveloping algebra of the corresponding Lie-Rinehart algebra.

In the supercase, the notion of Poisson superalgebra appeared in connection with supersymmetric mechanics. Indeed, the generalization of Hamiltonian mechanics to systems admitting ``bosonic" and ``fermionic" variables naturally led to the notion of Poisson superalgebra, see e.g. \cite{Kan,Kos, KS1, S, Sh}. Note that there is also a notion of \textit{Buttin bracket} (sometimes called \textit{antibracket} or \textit{Schouten bracket}), which is an odd version of the Poisson bracket, see e.g \cite{KS2,LS,S}. In a recent paper (\cite{La}), Lamkin studied the universal enveloping algebra of Poisson superalgebras. In particular, he constructed the universal enveloping algebra for a Lie-Rinehart superalgebra and proved a new version of the PBW Theorem, in order to generalize Huebschmann's constructions (see \cite{Hu}) to the supercase. Notably, he showed that the Poisson universal enveloping algebra can be constructed using the corresponding Lie-Rinehart superalgebra of K\"{a}hler differentials and proved a PBW Theorem for Poisson superalgebras.

\subsection{The case of the characteristic $p=2$} In the case where the base field has characteristic $p=2$, the usual definition of a Lie superalgebra has to be amended. Indeed, a Lie superalgebra is a $\Z/2\Z$-graded vector superspace $L=L_\ev\oplus L_\od$ such that the even part $L_\ev$ is a Lie algebra, the odd part is a two-sided $L_\od$-module due to symmetry, and a map $s:L_\od\rightarrow L_\ev$, called \textit{squaring} that satisfies some compatibility conditions, including a modified Jacobi identity; see Section \ref{sec:sLie-2}. The definition is due to Lebedev (\cite{LeD}), although earlier versions can be found in \cite{BMPZ,NR}. Due to the presence of the squaring, new methods are required to study those objects, and they have received some proper attention over the last few years. For instance, double extensions and Manin triples have been studied in \cite{BB1,BB2}, while the classification was investigated in \cite{BGL1,BGLLS,BLLS}. Derivations and deformations were studied in \cite{BGL2,BGLL}, while a full cohomology theory was introduced in \cite{BM} and used in \cite{BBE} to study Lagrangian extensions and left-symmetric structures.

The analogues of several types of superalgebras, such as orthogonal, Hamiltonian, Poisson, and contact superalgebras were studied by Lebedev in \cite{LeD}. He showed that there are four equivalence classes of non-degenerate symmetric bilinear forms on a linear superspace, denoted $I$, $I\Pi$, $\Pi I$, and $\Pi\Pi$. This gives rise to several Poisson brackets (resp. Buttin brackets) depending on whether the form is even (resp. odd); see \cite{LeD} for further details.\footnote{Note that two of the Poisson brackets in \cite{LeD} (namely, types $II$ and $I\Pi$) do \emph{not} define a Lie superalgebra.}

More recently, an algebraic definition of a Poisson superalgebra in characteristic $2$ was introduced by Petrogradsky and Shestakov in \cite[Section~2.4]{PS}. Their definition includes an additional compatibility condition between the squaring map and the associative supercommutative product, see Section \ref{def:poisson}.
We also mention that \textit{restricted} Poisson algebras (non-graded) were studied in \cite{BK,BYZ} in the characteristic $p\geq3$ case and in \cite{BEL} in the characteristic $p=2$ case. In this last paper, we notably adapted some of Huebschmann's constructions; namely, we built a cohomology for restricted Poisson algebras and restricted Lie-Rinehart algebras and showed that under a freeness hypothesis, they are isomorphic when considering the restricted Lie-Rinehart algebra of K\"{a}hler differentials. However, we will \textit{not} consider any restricted Lie algebras in the present paper.

\subsection{Outline of the paper and main results} The aim of this paper is to investigate the notion of Poisson superalgebras over a field of characteristic $p = 2$ and to study analogues of classical constructions such as representations, cohomology, formal deformations and universal enveloping algebras. In Section~\ref{sec:2}, we review the basics of Lie superalgebras in characteristic $2$ and their cohomology. In Section~\ref{sec:3}, we introduce the notion of a Lie-Rinehart superalgebra in characteristic $2$ (see Definition~\ref{sec:superlr2}), construct its universal enveloping algebra, and show that it satisfies an appropriate universal property (Proposition~\ref{prop:univ-LR}).
The core of the paper is contained in Section~\ref{sec:4}. We recall the definition of a Poisson superalgebra in \ref{def:poisson} and provide an example. We also show that the Poisson brackets defined by Lebedev (\cite{LeD}) indeed give rise to Poisson superalgebras. In addition, we investigate representations and construct the semidirect product of a Poisson superalgebra with a representation (Proposition~\ref{prop:semi-direct-poisson}). Furthermore, we build a full cohomology complex for Poisson superalgebras in Section~\ref{sec:cohopoiss} and consider formal deformations. In particular, we prove that infinitesimal deformations of a Poisson superalgebra are captured by the even part of the second cohomology group (Proposition~\ref{prop:equiv}) and we compute examples.
We define the universal enveloping algebra of a Poisson superalgebra and, analogously to \cite{Hu,BYZ,La,BEL}, we show that the module of K\"{a}hler differentials of a Poisson superalgebra has a Lie-Rinehart superalgebra structure (Proposition~\ref{prop:kahler-LR}) and that its universal enveloping algebra can be seen as the Poisson universal enveloping algebra (Theorem~\ref{thm:univ-iso}). Finally, we initiate the study of pre-Poisson superalgebras in Appendix~\ref{sec:5}.

Throughout the text, $\K$ denotes a field of characteristic $p = 2$, and we use the notation $\Z_2 = \Z/2\Z$ for the group of integers modulo $2$.

\section{Background}\label{sec:2}
\subsection{Lie superalgebras in characteristic $2$}\label{sec:sLie-2}

Let $V$ and $W$ be vector spaces. A map $s:V\rightarrow W$ is called a \textit{squaring} if 
\begin{equation}\label{eq:sq}
    s(\lambda x)=\lambda^2s(x),~\forall \lambda\in\K,~\forall x\in V,\end{equation}
and if the map
\[
V\times V \rightarrow W \qquad (x,y)\mapsto s(x+y)-s(x)-s(y)
\] is bilinear. 

Following \cite{BMPZ, LeD}\footnote{This definition however appeared way earlier in a paper of Nijenhuis and Richardson (\cite{NR}) under the name ``strongly anticommutative graded algebras'', but was not used (to our knowledge) until the aforementioned works.}, a \emph{Lie superalgebra} over a field $\K$ of characteristic $p=2$ is a $\Z_2$-graded vector space $L=L_{\ev}\oplus L_{\od}$ such that the even part $L_{\ev}$ is a Lie algebra, the odd part $L_{\od}$ is a $L_{\ev}$-module made two-sided by symmetry (the bracket on $L_{\ev}$ as well as the action of $L_{\ev}$ on $L_{\od}$ are given by the same symbol $[-, -]$), and a squaring  $s:L_{\od}\rightarrow L_{\ev}$, such that the bracket of two odd elements is given by:
    \begin{equation}\label{eq:billy}
        [x,y]:=s(x+y)-s(x)-s(y),~~ \forall x,y\in L_{\od}.
    \end{equation}
The Jacobi identity involving the squaring reads  
    \begin{equation}\label{jacobi2} 
        [s(x),y]=[x,[x,y]],~\forall x\in L_{\od},~\forall y\in L.    \end{equation}
Such a Lie superalgebra will be denoted by $\left(L,[-,-],s\right)$ or by $L$ if no confusion is possible.

Note that in this paper, when referring to a squaring on a $\Z_2$-graded Lie algebra $L$, we mean a map $L_\od\rightarrow L_\ev$ satisfying \eqref{eq:sq} such that \eqref{jacobi2} is satisfied.

\sssbegin{Example}\label{ex:asso} Let $A=A_{\ev}\oplus A_{\od}$ be an associative superalgebra in characteristic $p=2$. It can be endowed with a Lie superalgebra structure with the standard commutator and the squaring $s(a):=a^2,~\forall a\in A_{\od}$.   This applies in particular for $A=\text{End}(V)$, for $V$ any $\Z_2$-graded vector space.
\end{Example}

    Let $\bigl(L,[-,-],s\bigl)$ and $\bigl(L',[-,-]',s'\bigl)$ be two Lie superalgebras. An even linear map\\ $\varphi: L\rightarrow L'$ is called a \emph{Lie superalgebras morphism} if
   \begin{equation} \begin{array}{rcll}
        \varphi\left([x,y]\right)&=&[\varphi(x),\varphi(y)]',&~\forall x\in L_{\ev},~\forall y\in L,\\
        ~\varphi\circ s(x)&=&s'\circ\varphi(x),&~\forall x\in L_{\od}. 
    \end{array}\end{equation}

A linear map $d: L\rightarrow L$ is called a \emph{derivation} of $L$ if
   \begin{equation}\label{eq:deriv} \begin{array}{rcll}
         d([x,y])&=&[d(x),y]+[x,d(y)],&~\forall x\in L_{\ev},~y\in L,\\
        d(s(x))&=&[d(x),x],&~\forall y\in L_{\od}.  
    \end{array}\end{equation}
Note that the space $\Der(L)$ of derivations of a Lie superalgebra $L$ is itself a Lie superalgebra with the commutator bracket and the squaring $d\mapsto d^2.$

  A \emph{representation} $(\gamma,V)$ of a Lie superalgebra $\left(L,[-,-],s\right)$ is a $\Z_2$-graded vector space $V$ and an even Lie superalgebras morphism $\gamma: L\rightarrow \text{End}(L)$. In that case, $V$ is called a \emph{$L$-module}.

\subsubsection{Subalgebras and ideals}\label{sec:Lie-ideal}
    Let $\left(L,[-,-],s\right)$ be a Lie superalgebra and $H\subseteq  L$ be an homogeneous linear subspace.
    \begin{itemize}
        \item[($i$)] The subspace $H$ is a \emph{Lie subalgebra} if it is closed under the bracket and the squaring.
        \item[($ii$)] The subspace $H$ is an \emph{ideal} if it is closed under the squaring and if $[h,x]\in H,~\forall h\in H,~ x\in L.$
    \end{itemize}

\sssbegin{Lemma}\label{lem:centralsquare}
Let $(L,[-,-],s)$ be a Lie superalgebra. A map $\bar{s}:L_\od\rightarrow L_\ev$ is a squaring on $L$ satisfying the Jacobi identity if and only if there exists a 2-semilinear map $g:L_\od\rightarrow Z(L)_\ev$ such that $\bar{s}(x)=s(x)+g(x),~\forall x\in L_\od$.
\end{Lemma}
\begin{proof}
    Suppose that $\bar{s}=s+g$ is a squaring satisfying the Jacobi identity. Then, for all $x,y\in L_\od,$ and all $\lambda\in \K,$ we have
    $$g(\lambda x+y)=\lambda^2\bar{s}(x)+\bar{s}(y)+\lambda^2s(x)+s(y)+2[\lambda x,y]=\lambda^2g(x)+g(y).$$ Thus, the map $g$ is 2-semilinear. Moreover, since for all $z\in L$, we have 
    $$[g(x),z]=[\bar{s}(x)+s(x),z]=2[x,[x,z]]=0,$$ it follows that $g$ maps any $x\in L_\od$ into the center of $L$. Next, suppose that $g:L_\od\rightarrow Z(L)_\ev$ is 2-semilinear and let $\bar{s}=s+g.$ Clearly, we have $\bar{s}(\lambda x)=\lambda^2\bar{s}(x),~\forall x\in L_\od,~\forall \lambda\in \K.$ Moreover, for all $x,y\in L_\od$ and all $z\in L$, we have
    $$\bar{s}(x+y)=g(x)+g(y)+s(x)+s(y)+[x,y]=\bar{s}(x)+\bar{s}(y)$$ as well as
    $$[\bar{s}(x),z]=[g(x)+s(x),z]=[s(x),z]=[x,[x,z]]. $$ Therefore, the map $\bar{s}$ is a squaring satisfying the Jacobi identity.
\end{proof}


\sssbegin{Proposition} \label{prop:jacobson}
    Let $L_\ev$ be a Lie algebra with a basis $(e_1,\cdots,e_n)$ and $L_\od$ be a $L_\ev$-module with a basis  $(x_1,\cdots x_m)$. Suppose that for all $1\leq j\leq m$, there exists $y_j\in L_\ev$ such that $\ad_{y_j}=\ad^2_{x_j}.$ Then, $L=L_\ev\oplus L_\od$ is a Lie superalgebra.
\end{Proposition}

\begin{proof}
Let $x_1,\cdots x_m$ is a basis of $L_\od$ and define $s(x_j):=y_j~\forall~1\leq j\leq n$. Then, for any $\lambda_i\in\K$, we define
\begin{equation}\label{eq:ext-sq}
    s\Bigl(\sum_{i=1}^m\lambda_ix_i\Bigl)=\sum_{i=1}^m\lambda_i^2s(x_i)+\sum_{1\leq i<j\leq m}\lambda_i\lambda_j[x_i,x_j].
\end{equation}
Let us prove that the map $s$ satisfies the Jacobi identity. For all $y\in L$, we have
\begin{align*}
    \sum_{i\neq j}\lambda_i\lambda_j[x_i,[x_j,y]]&=\sum_{i< j}\lambda_i\lambda_j[x_i,[x_j,y]]+\sum_{i> j}\lambda_i\lambda_j[x_i,[x_j,y]]\\
    &=\sum_{i< j}\lambda_i\lambda_j[x_i,[x_j,y]]+\sum_{i< j}\lambda_i\lambda_j[x_j,[x_i,y]]=\sum_{i<j}\lambda_i\lambda_j[[x_i,x_j],y].
    \end{align*} Thus, it follows that
\begin{align*}
&\Bigl[s\bigl(\sum_{i=1}^m\lambda_ix_i\bigl),y\Bigl]+\Bigl[\sum_{i=1}^m\lambda_ix_i,\Bigl[\sum_{j=1}^m\lambda_jx_j,y\Bigl]\Bigl]\\
    =&~\cancel{\sum_{i}\lambda_i^2[s(x_i),y]}+\sum_{i<j}\lambda_i\lambda_j[[x_i,x_j],y]+\cancel{\sum_{i}\lambda_i^2[x_i,[x_i,y]]}+\sum_{i\neq j}\lambda_i\lambda_j[x_i,[x_j,y]]=0.
\end{align*}
Therefore, the map $s$ satisfies the Jacobi identity. The conclusion follows.
\end{proof}

\subsection{Cohomology of Lie superalgebras in characteristic 2}\label{sec:cohoLie} In this Section, we recall the construction of the cohomology for Lie superalgebras in characteristic 2 that was given in \cite{BM}. Earlier instances can be found in \cite{BGL2,BGLL}.

Let $\bigl(L,[-,-],s\bigl)$ be a Lie superalgebra and let $M$ be a $L$-module. We set $XC^0(L;M):=C_{\mathrm{CE}}^0 (L;M)$ and $XC^1(L;M):=C_{\mathrm{CE}}^1 (L;M)$\footnote{The subscript CE refers to the usual Chevalley-Eilenberg cohomology, see \cite[Section 6.1]{BM}}.

Let $n\geq2$, $\varphi\in C^n_{\ce} (L;M)$, $\omega:L_\od\otimes \wedge^{n-2} L\rightarrow M$, $\lambda\in\K$ and $x,z_2,\cdots,z_{n-1}\in L$. The pair $(\varphi,\omega)$ is a $n$-cochain if
		\begin{align}
			\omega(\lambda x, z_2,\cdots,z_{n-1})=&~\lambda^2\omega(x,z_2,\cdots,z_{n-1})\\
			\omega(x+y,z_2,\cdots,z_{n-1})=&~\omega(x,z_2,\cdots,z_{n-1})+\omega(y,z_2,\cdots,z_{n-1})\\\nonumber&~+\varphi(x,y,z_2,\cdots,z_{n-1}),\\
        (z_2,\cdots, z_{n-1})\mapsto &~\w(-,z_2,\cdots, z_{n-1}) \text{ is linear}.
        \end{align}
        
		We denote by $XC^n(L;M)$ the space of $n$-cochains of $L$ with values in $M$.

\subsubsection{Coboundary operators}\label{sec:diff.op.p=2}
Let $n\geq 2$. The coboundary maps
	$d^n:XC^n(L;M)\rightarrow XC^{n+1}(L;M)$ are given by $d^n(\varphi,\omega)=\bigl(d^n_{\mathrm{CE}}(\varphi),\delta^n(\omega)\bigl),$ where
	\begin{align*}	\delta^n\omega(x,z_2,\cdots,z_n):=&~x\cdot\varphi(x,z_2,\cdots,z_n)+\sum_{i=2}^{n}z_i\cdot\omega(x,z_2,\cdots,\hat{z_i},\cdots,z_n)\\	&+\varphi\bigl(s(x),z_2,\cdots,z_n\bigl)+\sum_{i=2}^{n}\varphi\bigl([x,z_i],x,z_2,\cdots,\hat{z_i},\cdots,z_n \bigl)\\
		&+\sum_{2\leq i<j\leq n}\omega\bigl(x,[z_i,z_j],z_2,\cdots,\hat{z_i},\cdots,\hat{z_j},\cdots,z_n  \bigl).
	\end{align*}
    For $n=0,1$, we define $d^0=d^0_{\mathrm{CE}}$ and 
	\begin{align*}
		d^1:XC^1(L;M)&\rightarrow XC^2(L;M)\\
		\varphi&\mapsto\bigl(d_{\mathrm{CE}}^1\varphi,\delta^1\varphi\bigl),~\text{where }\delta^1\varphi(x):=\varphi\bigl( s(x)\bigl)+x\cdot\varphi(x),~\forall x\in L.
	\end{align*}

The complex $\left(XC^n(L;M),d^n \right)_{n\geq 0}$ is a cochain complex, see \cite[Theorem 6.1]{BM}. The $n^{th}$ cohomology group of the Lie superalgebra $L$ in characteristic 2 is defined by	
		$$ XH^n(L;M):=XZ^n(L;M)/XB^n(L;M),$$ 
		with $XZ^n(L;M)=\Ker(d^n)$ the restricted $n$-cocycles and $XB^n(L;M)=\text{Im}(d^{n-1})$ the restricted $n$-coboundaries. 

A wide variety of examples of computations of this cohomology can be found in \cite{BGLL,BGL2,BM,BBE}.

\subsubsection{Hochschild 2-cocycles vs. Lie superalgebras 2-cocycles} Let $A$ be an associative superalgebra. Let $\mu:A\times A\rightarrow A$ be an \textit{Hochschild 2-cocycle}, that is, a bilinear map satisfying
\begin{equation}
    x\mu(y,z)+\mu(xy,z)+\mu(x,yz)+\mu(x,y)z,\quad\forall x,y,z\in A.
\end{equation} Consider the maps
    \begin{align*}
        \varphi_{\mu}(x,y)&=\mu(x,y)+\mu(y,x), &\forall x,y\in A;\\
        \w_{\mu}(x)&=\mu(x,x), &\forall x\in A_\od.
    \end{align*} Then, $(\varphi_{\mu},\w_{\mu})\in XZ^2(A;A)$, where $A$ is seen as a Lie superalgebra as in Example \ref{ex:asso}.

\section{Lie-Rinehart superalgebras in characteristic $2$}\label{sec:3}
In this section we introduce analogs of Lie-Rinehart superalgebras in characteristic $p=2$. We give the main definitions, construct the universal enveloping algebra and discuss its universal property.

A superalgebra $(A,\cdot)$ is called \textit{supercommutative} if
$$ ab=ba,~\forall a,b\in A; \quad \text{and}\quad a^2=0,~\forall a\in A_\od.$$
\subsection{Lie-Rinehart superalgebras}
\label{sec:superlr2}
    Let $A$ be an associative supercommutative superalgebra, $(L,[-,-],s)$ be a Lie superalgebra and $\rho:L\rightarrow \Der(A)$ be a morphism of Lie superalgebras. The triple $(A,L,\rho)$ is called a \emph{Lie-Rinehart superalgebra in characteristic $2$} if 
    \begin{enumerate}
        \item[$(i)$] $L$ is an $A$-module;
        \item[$(ii)$] the map $\rho$ is $A$-linear;
        \item[$(iii)$] the three following compatibility conditions are satisfied:
\begin{equation}\label{eq:SLR-cond}
\begin{cases}\begin{array}{rcll}
[x,ay]&=&a[x,y]+\rho_x(a)y,&~\forall a\in A,~\forall x,y\in L;\\
                    s(ax)&=&a^2s(x)+\rho_{ax}(a)x,&~\forall a\in A_\ev, ~\forall x\in L_\od;\\
                    s(ax)&=&\rho_{ax}(a)x,&~\forall a\in A_\od, \forall x\in L_\ev.
\end{array}\end{cases}
  \end{equation}      
  \end{enumerate}
For the sake of readability, the notations $\rho_x$ and $\rho(x)$ will be used interchangeably throughout the text.
\sssbegin{Lemma}\label{lem:ext-rinhart}
    Suppose that conditions \eqref{eq:SLR-cond} hold for $x$ in the basis of $L$. Then, conditions \eqref{eq:SLR-cond} hold for any $x\in L$.
\end{Lemma}
\begin{proof}\renewcommand{\qedsymbol}{}
    We prove the Lemma in the case where $a$ is even, $x$ is odd. Denote $x=\sum_i\lambda_ix_i,$ for $\lambda_i\in\K$ and $x_1,\cdots, x_m$ a basis of $L_\od$. We have
    \begin{align*}
&s\Bigl(a\sum_i\lambda_ix_i\Bigl)+a^2\sum_i\lambda_is(x_i)+\rho\Bigl(a\sum_i\lambda_ix_i\Bigl)(a)\sum_j\lambda_jx_j\\
        =&~\cancel{\sum_i\lambda_i^2s(ax_i)}+\sum_{i<j}\lambda_i\lambda_j[ax_i,ax_j]+\cancel{a^2\sum_i\lambda_i^2s(x_i)}+a^2\sum_{i<j}\lambda_i\lambda_j[x_i,x_j]\\
        &+~\cancel{\sum_{i}\lambda_i^2\rho(ax_i)(a)x_i}+\sum_{i<j}\rho(ax_i)(a)x_j+\sum_{j<i}\rho(ax_i)(a)x_j\\
        =&~\sum_{i<j}\lambda_i\lambda_j\bigl(a[ax_i,x_j]+\rho(ax_i)(a)x_j\bigl)+a^2\sum_{i<j}\lambda_i\lambda_j[x_i,x_j]\\
        &+~\sum_{i<j}\lambda_i\lambda_j
        \rho(ax_i)(a)x_j+\sum_{i<j}\lambda_i\lambda_j
        \rho(ax_j)(a)x_i \\
        =&~\sum_{i<j}\lambda_i\lambda_j\bigl(a^2[x_i,x_j]+a\rho(x_j)(a)x_i+a\rho(x_i)(a)x_j\bigl)\\
        &+~a^2\sum_{i<j}\lambda_i\lambda_j[x_i,x_j]+\sum_{i<j}\lambda_i\lambda_j\rho(ax_i)(a)x_j+\sum_{i<j}\lambda_i\lambda_j\rho(ax_j)(a)x_i=0.\quad\square\\
\end{align*}
\end{proof}

\sssbegin{Example}[Lie superalgebras of derivations] Let $A$ be an associative supercommutative superalgebra and let $\Der(A)$ be its Lie superalgebra of derivations (see Eq. \ref{eq:deriv}). Then, $\Der(A)$ is an $A$-module with the action given by $(a\cdot D)(b):=aD(b),~\forall a,b\in A,~\forall D\in \Der(A).$ As anchor map, we take the identity map $\id:\Der(A)\rightarrow \Der(A).$ It is indeed an $A$-linear superalgebras morphism. The first compatibility condition of Definition \ref{sec:superlr2} is readily checked (see e.g.  \cite{R}). Let $(a,b)\in A_\ev\times A$ and $D\in\Der(A)_\od$. We have
\small{}\begin{align*}
    s\bigl((a\cdot D)\bigl)(b)=(a\cdot D)\bigl(aD(b)\bigl)=aD(a)D(b)+a^2(D\circ D)(b)=a^2s(D)(b)+\bigl((a\cdot D)(a)\cdot D\bigl)(b).
\end{align*}\normalsize{}
In the case where $a\in A_\od$ and $D\in \Der(A)_\ev,$ the same computation leads to $s\bigl((a\cdot D)\bigl)(b)=\bigl((a\cdot D)(a)\cdot D\bigl)(b)$, since $a^2=0$. The second compatibility condition is therefore satisfied and $\bigl(A,\Der(A),\id\bigl)$ is a Lie-Rinehart superalgebra.
\end{Example}

 Let $(A,L,\rho)$ and $(A',L',\rho')$  be Lie-Rinehart superalgebras in characteristic $2$. A morphism of Lie-Rinehart superalgebras is a pair $(\phi,\psi)$, with $\phi:A\rightarrow A'$ is a morphism of associative superalgebras and $\psi:L\rightarrow L'$ is a $A$-linear Lie superalgebras morphism, such that
    $$\phi\left(\rho_x(a)\right)=\rho'_{\psi(x)}\left(\phi(a)\right),~~\forall x\in L, \forall a\in A.$$
    
    Let $(A,L,\rho)$ be a Lie-Rinehart superalgebra. An ideal $I$ of the Lie superalgebra $L$ is called \textit{ideal of the Lie-Rinehart superalgebra} $(A,L,\rho)$ if $(A,I,0)$ is a Lie-Rinehart superalgebra with the induced action and the anchor identically equal to zero. 
 With this definition, one can check that $\bigl(A,L/I,\Tilde{\rho}\bigl)$ is a Lie-Rinehart superalgebra with $\Tilde{\rho}$ the induced anchor.
\sssbegin{Example} Let $(A,L,\rho)$ be a Lie-Rinehart superalgebra. Then, the anchor map induces a Lie-Rinehart morphism $(\id,\rho):(A,L,\rho)\rightarrow\bigl(A,\Der(A),\id\bigl)$, whose kernel is an ideal.
\end{Example}

\subsubsection{Lie-Rinehart modules}
    Let $(A,L,\rho)$ be a Lie-Rinehart superalgebra and $V$ be a $A$-module. A \textit{weak module} is a pair $(\pi,V)$, where $\pi:L\rightarrow\text{End}(V)$ is a Lie superalgebras morphism and $V$ is a $\Z_2$-graded vector space such that
    \begin{equation}
                \pi_x(av)=a\pi_x(v)+\rho_x(a)v,~\forall a\in A,~\forall x\in L,~\forall v\in V. 
    \end{equation}
  If the map $\pi$ is $A$-linear, the pair $(\pi,V)$ is called \textit{strong module}.

  \sssbegin{Remark} The above definition resembles the classical one (see \cite{Ri,Hu}), but differs from \cite[Definition 5.4]{La}. Our definition is motivated by the following Proposition.
  \end{Remark}

\sssbegin{Proposition} Let $(A,L,\rho)$ be a Lie-Rinehart superalgebra and $V$ be a strong Lie-Rinehart module. Then, the triple $(A,L\oplus V,\Tilde{\rho})$ is a Lie-Rinehart superalgebra with
 \begin{itemize}
        \item[$(i)$] the bracket $\bigl[(x+v),(y+w)\bigl]:=[x,y]+\pi_x(w)+\pi_y(v),~\forall x,y\in L,~\forall v,w\in V;$
        \item[$(ii)$] the squaring $ s(x+v):=s(x)+\pi_x(v),~ \forall x\in L_\od,~ \forall v\in V_\od;$
        \item[$(iii)$] the anchor $\Tilde{\rho}(x+v)=\rho(x),~\forall x\in L~\forall v\in V.$
    \end{itemize}
\end{Proposition}
\begin{proof}
Straightforward computations.
\end{proof}

\subsection{Universal enveloping algebra}
This section is devoted to the construction of the universal enveloping algebra of a Lie-Rinehart superalgebra. Our construction follows the original work of Rinehart (\cite{Ri}) and \cite[Section 5]{La}. The universal property for Lie-Rinehart algebras was proven later by Huebschmann (\cite{Hu}). Given a Lie-Rinehart superalgebra $(A,L,\rho)$, the idea is to equip the superspace $A \oplus L$ with a Lie superalgebra structure and then to consider its enveloping algebra.

\sssbegin{Proposition}\label{prop:direct-sum}
    Let $(A,L,\rho)$ be a Lie-Rinehart superalgebra. Then, the superspace $A\oplus L$ is a Lie superalgebra with 
    \begin{itemize}
        \item[$(i)$] the bracket $\bigl[a+x,b+y\bigl]:=\rho_x(b)+\rho_y(a)+ [x,y] ,~\forall a,b\in A,~\forall x,y\in L;$
        \item[$(ii)$] the squaring $ s(a+x):=\rho_x(a)+s(x),~ \forall a\in A,~ x\in L.$
    \end{itemize}
    
\end{Proposition}

\begin{proof}
    Recall that the grading on $A\oplus L$ is given by $(A\oplus L)_{\ev}=A_{\ev}\oplus L_{\ev}$ and $(A\oplus L)_{\od}=A_{\od}\oplus L_{\od}$. We will only prove the condition on the squaring (for the bracket, see \cite{Ri} or \cite{La}). Let $a+x\in A_{\od}\oplus L_{\od}$ and $\lambda\in \K$. We have
    $$s\bigl(\lambda(a+x)\bigl)=s\bigl(\lambda a+\lambda x\bigl)=\bigl(\rho_{\lambda x}(\lambda a)+s(\lambda x)\bigl)=\lambda^2\bigl(\rho_x(a)+s(x)\bigl)=\lambda^2s(a+x).$$
Moreover, for any $b+y\in A_{\od}\oplus L_{\od}$, we have  
    \begin{align*}
        s\bigl((a+x)+(b+y)\bigl)&=s(a+b+x+y)\\
        &=\bigl(\rho_{a+b}(x+y)+ s(x+y)   \bigl)\\
        &=\bigl(\rho_a(x)+s(x)\bigl)+\bigl(\rho_b(y)+s(y)\bigl)+ \bigl(\rho_a(y)+\rho_b(x)+[x,y]\bigl)\\
        &=s(a+x)+s(b+y)+\bigl[a+x,b+y\bigl].
    \end{align*}
    It remains to check the Jacobi identity involving the squaring (see Eq. (\ref{jacobi2})). Let $a+x\in A_{\od}\oplus L_{\od}$ and $b+y\in A\oplus L$.
    On one hand, we have
    \begin{align*}
        \bigl[s(a+x),b+y\bigl]&=\bigl[\rho_x(a)+s(x),b+y \bigl]\\
        &=\bigl(\rho_{s(x)}(b)+\rho_y\bigl(\rho_x(a)\bigl)+\bigl[s(x),y\bigl] \bigl)\\
        &=\bigl(\rho_x^2(b)+\rho_y\circ\rho_x(a)+[x,[x,y]]\bigl).
    \end{align*}
    On the other hand, we have
    \begin{align*}
        \Bigl[a+x,\bigl[a+x,b+y\bigl]\Bigl]&=\Bigl[(a+x),\bigl(\rho_x(b)+\rho_y(a)+[x,y]\bigl)\Bigl]\\
        &=\bigl(\rho_x\bigl(\rho_x(b)+\rho_y(a)\bigl)+\rho_{[x,y]}(a)+\bigl[x,[x,y]\bigl]\bigl)\\
        &=\bigl(\rho_x^2(b)+\underset{=~0}{\underbrace{\rho_x\circ\rho_y(a)+\rho_x\circ\rho_y(a)}}+\rho_y\circ\rho_x(a)+\bigl[x,[x,y]\bigl]\bigl)\\
        &=\bigl(\rho_x^2(b)+\rho_y\circ\rho_x(a)+\bigl[x,[x,y]\bigl]\bigl).
    \end{align*}
    Therefore, we have $$\bigl[s(a+x),b+y\bigl]=\bigl[a+x,\bigl[a+x,b+y\bigl]\bigl].$$
    It follows that Eq. (\ref{jacobi2}) is satisfied and $A\oplus L$ is a Lie superalgebra.
\end{proof}

\subsubsection{Construction of the universal enveloping algebra.} Let $(A,L,\rho)$ be a Lie-Rinehart superalgebra. Denote by $T(A\oplus L)=\displaystyle\bigoplus_{n\geq 0}(A\oplus L)^{\otimes n}$ the tensor algebra of $A\oplus L$. 
Let $U(A\oplus L)=\bigl(T(A\oplus L)/I_1\bigl)/I_2$ be the Lie enveloping algebra of $A\oplus L$, with $I_1$ the two-sided ideal of $T(A\oplus L)$ generated by elements
$$I_1:=\Bigl<(a+x)\otimes(b+y)+(b+y)\otimes(a+x)+\bigl[(a+x),(b+y)\bigl],~(a+x),~(b+y)\in (A\oplus L)\Bigl>$$ and $I_2$ the two-sided ideal of $T(A\oplus L)$ generated by elements
$$I_2:=\Bigl<s(c+z)+(c+z)\otimes(c+z),~(c+z)\in (A\oplus L)_\od\Bigl>,$$
see also \cite[Section 2.2]{PS}. Let $i:A\oplus L\rightarrow U(A\oplus L)$ be the canonical inclusion, let $i_A: A\rightarrow A\oplus L$ and $i_L: L\rightarrow A\oplus L$ be the inclusions of $A$ and $L$ into $A\oplus L$ respectively. Moreover, let $\overline{U}(A\oplus L)$ be the subalgebra of $U(A\oplus L)$ generated by $i(A\oplus L)$. We define the\textit{ universal enveloping algebra} of $(A,L,\rho)$, denoted by $U(A,L,\rho)$ by
$$U(A,L,\rho)=\overline{U}(A\oplus L)/J,$$ with $J$ the two-sided ideal generated by elements
 $$J:=\bigl<i(a)i(b+x)+i(ab+ax)\bigl>.$$ 
Note that the following relations hold in $U(A,L,\rho)$ for all $a\in A$ and all $x\in L$:
 \begin{align}
     i_L(ax)&=i_A(a)i_L(x);\\
     i_L(x)i_A(a)&=i_A(a)i_L(x)+i_A\bigl(\rho_x(a)\bigl).
 \end{align}

The enveloping algebra $U(A,L,\rho)$ satisfies the following universal property.

\sssbegin{Proposition}\label{prop:univ-LR}
For any triple $(B,r_A,r_L)$ where $B$ is an associative superalgebra,
$r_A:A\to B$ is an even superalgebra morphism, and
$r_L:L\to B$ is an even Lie superalgebra morphism satisfying for all $a\in A$ and all $x\in L$
 \begin{align}
     r_L(ax)&=r_A(a)r_L(x);\\
     r_L(x)r_A(a)&=r_A(a)r_L(x)+r_A\bigl(\rho_x(a)\bigl),
 \end{align} there exists an unique even associative
superalgebras morphism $f:U(A,L,\rho)\to B$ such that $f\circ i_A=r_A$ and $f\circ i_L=r_L$.
\end{Proposition}

\sssbegin{proof} Consider the map $r:A\oplus L\to B,~a+x\mapsto r_A(a)+r_L(x).$ It follows from \cite[Proposition 5.5]{La} that
$$ r\bigl([a+x,b+y]\bigl)=r(a+x)r(b+y)+r(b+y)r(a+x).$$ Moreover, for all $a+x\in (A\oplus L)_\od$, we have
\begin{align*}
    r(a+x)^2&=r_A(a)^2+r_A(a)r_L(x)+r_L(x)r_A(a)+r_L(x)^2\\
            &=r_A\bigl(\rho_x(a)\bigl)+r_L\bigl(s(x)\bigl)\\
            &=r\bigl(\rho_x(a)+s(x)\bigl)=r\bigl(s(a+x)\bigl).
\end{align*} Therefore, the map $r$ is a morphism of Lie superalgebras. By the universal property of $U(A\oplus L)$, there exists an unique associative superalgebras morphism $\widetilde{r}:\overline{U}(A\oplus L)\to B$ such that $\widetilde{r}\circ i=r$, where $i:A\oplus L\to\overline{U}(A\oplus L)$ denotes the canonical map.
Let $a,b\in A$ and $x\in L$. We have
\begin{align*}
    \widetilde{r}\bigl(i(a)i(b,x)+i(ab,ax)\bigl)&=\widetilde{r}(i(a))\widetilde{r}(i(b,x))+r(i(ab,ax))\\
    &=r_A(a)\bigl(r_A(b)+r_L(x)\bigl)+r_A(ab)+r_L(ax)=0,
\end{align*} since $r_A$ is a morphism of associative superalgebras and $r_L(ax)=r_A(a)r_L(x).$ Thus, the ideal $J$ is contained in the kernel of the map $\widetilde{r}$ and it follows that there exists an unique superalgebras morphism $f:U(A,L,\rho)\to B$ that satisfies $f\circ i_A=r_A$ and $f\circ i_L=r_L$.
\end{proof}

\section{Poisson superalgebras in characteristic $2$}\label{sec:4}

In this Section, we investigate representations, cohomology, formal deformations, universal enveloping algebra and examples of Poisson superalgebras.

\subsection{Poisson superalgebras}\label{def:poisson} 
    Let $\bigl(P,\{-,-\},s\bigl)$ be a Lie superalgebra equipped with an associative supercommutative product denoted by juxtaposition. Following \cite[Section 2.4]{PS}, $\bigl(P,\{-,-\},s\bigl)$ is called \textit{Poisson superalgebra in characteristic 2} if
    $$\begin{array}{crcll}
       (i)&\{xy,z\}&=&x\{y,z\}+\{x,z\}y,&~\forall x,y,z\in P;\\
        (ii)& s(xy)&=&x^2s(y)+xy\{x,y\},&~\forall x\in P_\ev,~\forall y\in P_\od.
    \end{array}$$
A Poisson superalgebra will be usually denoted by a triple $\bigl(P,\{-,-\},s\bigl)$, or sometimes by $\bigl(P,\cdot,\{-,-\},s\bigl)$ if it is needed to specify the associative supercommutative product $(\cdot)$. In the case where no confusion is possible, we denote it simply by $P$.\\
Let $P$ and $Q$ be Poisson superalgebras. A linear map $\varphi:P\rightarrow Q$ is called \textit{Poisson superalgebras morphism} if it is a morphism of associative superalgebras and of Lie superalgebras at the same time.

 Let $(P,\{-,-\},s)$ be a Poisson superalgebra. A linear map $\delta: P\to P$ is called a \textit{Poisson derivation} if 
$$\begin{array}{crcll}
    (i)& \delta(xy)&=&\delta(x)y+x\delta(y),&\forall x,y\in P,\\
    (ii)& \delta\bigl(\{x,y\}\bigl)&=&\{\delta(x),y\}+\{x,\delta(y)\},&\forall x,y\in P,\\
    (iii)& \delta\bigl(s(x)\bigl)&=&\{\delta(x),x\},&\forall x\in P_\od.
\end{array}$$ The space of Poisson derivations will be denoted by $\Der_{\mathrm{po}}(P).$ The space of linear maps $P\to P$ satisfying Condition $(i)$ only will be denoted by $\Der_{\mathrm{asso}}(P).$ From the definitions, it is clear that a Lie superalgebra derivation $\delta$ (see Eq. \eqref{eq:deriv}) is a Poisson derivation if and only if $\delta$ satisfies Condition ($i$) above.

An homogeneous Lie ideal (see  Section \ref{sec:Lie-ideal}) $I\subset P$ of a Poisson superalgebra is a \textit{Poisson ideal} if $xI\in I,~\forall x\in P.$ With this definition, one can readily check that the quotient $P/I$ is a Poisson superalgebra.

\sssbegin{Example}\label{ex:example} Consider the Lie superalgebra $P$ given by the basis $e_1,e_2|e_3,e_4$ (even$|$odd) and the squaring $s(e_3)=e_1,~s(e_4)=e_2,~s(e_3+\lambda e_4)=e_1+\lambda^2 e_2,$ for all $\lambda\in\K$. The brackets are zero. Consider the associative supercommutative product given by $e_1x=x,~\forall x\in P$, and $e_3e_4=e_4e_3=e_2$. Then, $P$ has the structure of a Poisson superalgebra in characteristic $2$.
\end{Example}

\sssbegin{Proposition}\label{tensorLie}
    Let $(P,\{-,-\},s)$ and $(Q,[-,-],s')$ be Poisson superalgebras. The tensor product $P\otimes Q$ is a Poisson superalgebra with
    \begin{enumerate}
           \item[($i$)] the associative supercommutative product $$(p_1\otimes q_1)(p_2\otimes q_2):=(p_1p_2)\otimes(q_1q_2),\quad \forall p_1,p_2\in P,~ \forall q_1,q_2\in Q; $$
           \item[($ii$)] the Lie bracket $$\bigl[p_1\otimes q_1,p_2\otimes q_2\bigl]_{\otimes}:=\{p_1,p_2\}\otimes q_1q_2+p_1p_2\otimes[q_1,q_2],\quad \quad \forall p_1,p_2\in P,~ \forall q_1,q_2\in Q;$$
        \item[($iii$)] the squaring 
        $$\begin{array}{rcll}
        s_{\otimes}(p\otimes q)&:=&p^2\otimes s'(q),&\forall p\in P_\ev,~\forall q\in Q_\od,\\
        s_{\otimes}(p\otimes q)&:=&s(p) \otimes q^2,&\forall p\in P_\od,~\forall q\in Q_\ev.
        \end{array}$$
    \end{enumerate}
\end{Proposition}
    
\begin{proof}
We prove the identities involving the squaring. Recall that the grading on $P\otimes Q$ is given by
$$(P\otimes Q)_\ev=(P_\ev\otimes Q_\ev)\oplus(P_\od\otimes Q_\od);~~(P\otimes Q)_\od=(P_\ev\otimes Q_\od)\oplus(P_\od\otimes Q_\ev).$$
Let $p_2\in P,~q_2\in Q$, $p_1\in P_\ev$ and $q_1\in Q_\od$. We have
\begin{align*}
    [s_{\otimes}(p_1\otimes q_1),p_2\otimes q_2]_{\otimes}&=[p_1^2\otimes s'(q_1),p_2\otimes q_2]_{\otimes}\\
    &=\{p_1^2,p_2\}\otimes s'(q_1)q_2+p_1^2p_2\otimes[s'(q_1),q_2]\\
    &=\bigl(p_1\{p_1,p_2\}+\{p_1,p_2\}p_1\bigl)\otimes s'(q_1)q_2+p_1^2p_2\otimes[q_1,[q_1,q_2]]\\
    &=p_1^2p_2\otimes[q_1,[q_1,q_2]].
\end{align*} On the other hand, we have
\begin{align*}
    \bigl[p_1\otimes q_1,[p_1\otimes q_1,p_2\otimes q_2]_{\otimes}\bigl]_{\otimes}=&~\bigl[p_1\otimes q_1,\{p_1,p_2\}\otimes q_1q_2+p_1p_2\otimes[q_1,q_2]\bigl]_{\otimes}\\
    =&~\{p_1,\{p_1,p_2\}\}\otimes q_1^2q_2+p_1\cancel{\{p_1,p_2\}\otimes[q_1,q_1q_2]}\\
    &+\cancel{\{p_1,p_1p_2\}\otimes q_1[q_1,q_2]}+p_1^2p_2\otimes[q_1,[q_1,q_2]]\\
    =&~p_1^2p_2\otimes [q_1,[q_1,q_2]].
\end{align*} Therefore, the squaring $s_{\otimes}$ satisfies the Jacobi identity and $P\otimes Q$ is a Lie superalgebra.
It remains to check the compatibility of the squaring $s_\otimes$ with the multiplication. There are two cases to consider.

\underline{The case where $p_2,q_2$ are even.} We have
\begin{align*}
    s_\otimes\bigl((p_1\otimes q_1)(p_2\otimes q_2)\bigl)=&~s_\otimes(p_1p_2\otimes q_1q_2)=(p_1p_2)^2\otimes s'(q_1q_2)\\
    =&~(p_1p_2)^2\otimes q_2^2s'(q_1)+q_1q_2[q_1,q_2]\\
    =&~(p_1p_2)^2\otimes q_2^2s'(q_1)+(p_1p_2)^2\otimes q_1q_2[q_1,q_2]\\
    =&~(p_2\otimes q_2)^2\bigl(p_1^2\otimes s'(q_1)\bigl)+(p_1p_2\otimes q_1q_2)\bigl(\{p_1,p_2\}\otimes q_1q_2+p_1p_2[q_1,q_2]\bigl)\\
    =&~(p_2\otimes q_2)^2s_{\otimes}(p_1\otimes q_1)+(p_1\otimes q_1)(p_2\otimes q_2)\bigl[p_1\otimes q_1,p_2\otimes q_2\bigl]_{\otimes}.
\end{align*}

\underline{The case where $p_2,q_2$ are odd.} In that case, we have
$$s_\otimes\bigl((p_1\otimes q_1)(p_2\otimes q_2)\bigl)=0=(p_2\otimes q_2)^2s_{\otimes}(p_1\otimes q_1)+(p_1\otimes q_1)(p_2\otimes q_2)\bigl[p_1\otimes q_1,p_2\otimes q_2\bigl]_{\otimes}.$$ It follows that $P\otimes Q$ is a Poisson superalgebra.
\end{proof}

\subsection{Analogs of ``classical" Poisson Lie superalgebras from \cite{LeD}}

In this Section, we show that the analogs of classical Poisson Lie superalgebras in characteristic 2 defined by Lebedev in \cite{LeD} are indeed Poisson superalgebras according to Definition \ref{def:poisson}.

 Let $k\geq 1$. The associative commutative algebra of divided powers in $k$ variables $x:=(x_1,\cdots,x_k)$ is defined for $\underline{N}:=(n_1,\cdots,n_k)$, where $n_s\geq 0,~\forall~  1\leq s \leq k$, by (see, e.g,  \cite{Ko, SF})
\begin{align}
    \K(x;\underline{N}):=\Span\bigl\{x^{(\underline{i})}:=x_1^{(i_1)}\cdots x_k^{(i_k)},~(\underline{i})=(i_1,\cdots,i_k),~0\leq i_s\leq 2^{n_s}-1  \bigl\}.
\end{align}
The multiplication is given by
\begin{equation}\label{eq:divpwr}
x^{(\underline{i})}x^{(\underline{j})}=\binom{\underline{i}+\underline{j}}{\underline{i}}x^{(\underline{i}+\underline{j})},~\text{where}~\binom{\underline{i}+\underline{j}}{\underline{i}}:=\prod_{s=1}^k\binom{i_s+j_s}{i_s}.
\end{equation}
We will write $\mathcal{O}(k,\underline{N})=\K(x;\underline{N})$ or just $\mathcal{O}(k)$ if the value of $\underline{N}$ is not important, see \cite{LeD} for more details. It is shown in \cite[Section 6]{LeD} that there are two Poisson Lie superalgebras in characteristic 2 given by the following formulas, for $f,g\in \mathcal{O}(n)$ and $h\in \mathcal{O}(n)_\od$.

\begin{enumerate}
    \item[$\textup{(}i\textup{)}$] \underline{Type $\Pi I$}. The bracket is given by
    \begin{equation}
        \{f,g\}_{\Pi I}=\sum_{i=1}^{k_\ev}\Bigl(\frac{\del f }{\del p_i}\frac{\del g }{\del q_i}+\frac{\del f }{\del q_i}\frac{\del g }{\del p_i}\Bigl)+\sum_{i=1}^{n_\od}\frac{\del f }{\del \theta_i}\frac{\del g }{\del \theta_i}
    \end{equation} and the squaring by
       \begin{equation}
        s_{\Pi I}(h)=\sum_{i=1}^{k_\ev}\Bigl(\frac{\del h }{\del p_i}\frac{\del h }{\del q_i}\Bigl)+\sum_{i=1}^{n_\od}\Bigl(\frac{\del h }{\del \theta_i}\Bigl)^{(2)};
    \end{equation}
    \item[$\textup{(}ii\textup{)}$] \underline{Type $\Pi\Pi$}. The bracket is given by \begin{equation}
        \{f,g\}_{\Pi \Pi}=\sum_{i=1}^{k_\ev}\Bigl(\frac{\del f }{\del p_i}\frac{\del g }{\del q_i}+\frac{\del f }{\del q_i}\frac{\del g }{\del p_i}\Bigl)+\sum_{i=1}^{k_\od}\Bigl(\frac{\del f }{\del \xi_i}\frac{\del g }{\del \eta_i}+\frac{\del f }{\del \eta_i}\frac{\del g }{\del \xi_i}\Bigl)
    \end{equation} and the squaring by
       \begin{equation}
        s_{\Pi \Pi}(h)=\sum_{i=1}^{k_\ev}\Bigl(\frac{\del h }{\del p_i}\frac{\del h }{\del q_i}\Bigl)+\sum_{i=1}^{k_\od}\frac{\del h }{\del \xi_i}\frac{\del h }{\del \eta_i};
        \end{equation} where we denoted $n=n_\ev+n_\od$, with $n_\ev$ (resp. $n_\od$) the dimension of the even part (resp. odd), and 
        \begin{equation*}
        \begin{array}{lll}
            p_i=x_i,&q_i=x_{k_\ev+i},&\text{if $n_\ev=2k_\ev$ and $1\leq i\leq k_\ev$;}\\
            \theta_i=x_{n_\ev+i},&&\text{if $1\leq i\leq n_\ev$;}  \\       
            \xi_i=x_{n_\ev+i},&\eta_i=x_{n_\ev+k_\od+i},&\text{if $n_\od=2k_\od$ and $1\leq i\leq k_\od$.}
        \end{array}\end{equation*}
\end{enumerate}

\sssbegin{Proposition}
    The Poisson Lie superalgebras of type $\Pi I$ and $\Pi\Pi$ are Poisson superalgebras in the sense of Definition \ref{def:poisson}.
\end{Proposition}
\begin{proof}
We prove that the Poisson Lie superalgebra of type $\Pi I$ satisfies Condition ($ii$) of Definition \ref{def:poisson}. The case of type $\Pi\Pi$ is similar. Let $f\in \mathcal{O}(n)$ and let $g\in \mathcal{O}(n)_\od$. We have (for the sake of readability, we drop the subscript $\Pi I$)
\begin{align*}
    s(fg)&=\sum_{i=1}^{k_\ev}\Bigl(g\frac{\del f}{\del p_i}+f\frac{\del g}{\del p_i}\Bigl)\Bigl(g\frac{\del f}{\del q_i}+f\frac{\del g}{\del q_i}\Bigl)+\sum_{i=1}^{n_\od}\Bigl(g\frac{\del f}{\del \theta_i}+f\frac{\del g}{\del\theta_i}\Bigl)\\
    &=\sum_{i=1}^{k_\ev}\Bigl(gf\frac{\del f}{\del p_i}\frac{\del g}{\del q_i}+fg\frac{\del g}{\del p_i}\frac{\del f}{\del q_i}+f^{(2)}\frac{\del g}{\del p_i}\frac{\del g}{\del q_i}\Bigl)+\sum_{i=1}^{n_\od}\Biggl(f^{(2)}\Bigl(\frac{\del g}{\del\theta_i}\Bigl)^{(2)}+fg\frac{\del f}{\del \theta_i}\frac{\del g}{\del \theta_i}\Biggl)\\
    &=f^{(2)}\Biggl(\sum_{i=1}^{k_\ev}\frac{\del g}{\del p_i}\frac{\del g}{\del q_i}+\sum_{i=1}^{n_\od}\Bigl(\frac{\del g}{\del\theta_i}\Bigl)^{(2)}\Biggl)+fg\Biggl(\sum_{i=1}^{k_\ev}\Bigl(\frac{\del f}{\del p_i}\frac{\del g}{\del q_i}+\frac{\del g}{\del p_i}\frac{\del f}{\del q_i}\Bigl)+\sum_{i=1}^{n_\od}\frac{\del f}{\del \theta_i}\frac{\del g}{\del \theta_i}\Biggl)\\
    &=f^{(2)}s(g)+fg\{f,g\}.
\end{align*}
\end{proof}
\subsection{Representations of Poisson superalgebras} Let $\bigl(P,\{-,-\},s\bigl)$ be a Poisson superalgebra and let $V=V_\ev\oplus V_\od$ be a super vector space. A \textit{representation} of $P$ in $V$ is given by two maps $\pi,\gamma:P\rightarrow\End(V)$ such that
\begin{enumerate}
    \item[($i$)] $\pi$ is a morphism of associative superalgebras;
    \item[($ii$)] $\gamma$ is a morphism of Lie superalgebras;
    \item[($iii$)] the four following compatibility conditions are satisfied:
    \begin{align}
    \gamma(xy)&=\pi(x)\gamma(y)+\pi(y)\gamma(x),&\forall x,y\in P;\label{eq:rep1}\\
    \pi\bigl([x,y]\bigl)&=\pi(x)\gamma(y)+\gamma(y)\pi(x),& \forall x,y\in P;\label{eq:rep2}\\
    \pi(y)\gamma(x)\pi(x)_{|V_\ev}&=\pi(x)\pi(y)\gamma(x)_{|V_\ev},&\forall x\in P_\od,~\forall y\in P_\ev;\label{eq:rep3}\\
    \pi(y)^2\gamma(x)_{|V_\od}&=\pi(y)\gamma(x)\pi(y)_{|V_\od}+\pi(x)\pi(y)\gamma(y)_{|V_\od}\label{eq:rep4}\\\nonumber&~~\quad+\gamma(y)\pi(x)\pi(y)_{|V_\od},& \forall x\in P_\od,~\forall y\in P_\ev.
    \end{align}
\end{enumerate}
The compatibility conditions \eqref{eq:rep1} and \eqref{eq:rep2} appeared in \cite{Oh, LOV, La}, while \eqref{eq:rep3} and \eqref{eq:rep4} are specific to the characteristic 2 case. Such a representation will be denoted by a triple $(V,\pi,\gamma)$.

\sssbegin{Example}[Adjoint representation]
    Let $\bigl(P,\{-,-\},s\bigl)$ be a Poisson superalgebra. Define $\pi(x)(v)=xv$ and $\gamma(x)(v)=\{x,v\},\quad\forall x,v\in P$. Then $(P,\pi,\gamma)$ is a representation of the Poisson superalgebra $P$ into itself, called \emph{adjoint representation}.
\begin{proof}
    We prove conditions \eqref{eq:rep3} and \eqref{eq:rep4}. Let $x,v\in P_\od$ and $y,w\in P_\ev$. Since $[x,x]=0$, \eqref{eq:rep3} is immediate. For \eqref{eq:rep4}, notice that
    \begin{equation}
        y^2\bigl(s(x+v)+s(x)+s(v)+[x,v]\bigl)=0=y\bigl(v[x,y]+v[x,y]\bigl).
    \end{equation} Thus, we have
    \begin{align*}       y^2[x,v]&=y^2s(x+v)+y^2s(x)+y^2s(v)\\&=y^2[x,v]+yv[x,y]+\underset{=~0}{\underbrace{xy[y,v]+xy[y,v]}}+\underset{=~vy[y,x]}{\underbrace{v[y,xy]}},
    \end{align*} which is exactly \eqref{eq:rep4}.
\end{proof}
\end{Example}

\sssbegin{Proposition}[Semi-direct product]\label{prop:semi-direct-poisson}
    Let $\bigl(P,\{-,-\},s\bigl)$ be a Poisson superalgebra and let $(V,\pi,\gamma)$ be a representation. Then, $\bigl(P\oplus V,\cdot,\{-,-\}_{\rtimes},s_{\rtimes}\bigl)$ is a Poisson superalgebra, where
    \begin{align}
        (x+v)\cdot (y+w)&=xy+\pi(x)(w)+\pi(y)(v),&\forall x,y\in P,\forall v,w\in V;\\
        \{x+v,y+w\}_{\rtimes}&=\{x,y\}+\gamma(x)(w)+\gamma(y)(v),&\forall x,y\in P,\forall v,w\in V;\\
        s_{\rtimes}(x+v)&=s(x)+\gamma(x)(v),&\forall x\in P_\od,\forall v\in V_\od.
    \end{align}
\end{Proposition}
\begin{proof}
    We prove that $s_{\rtimes}$ satisfy condition $(ii)$ of Definition \ref{def:poisson}. For the other conditions involving the squaring, see \cite[Theorem 3.1]{BM}. Let $x\in P_\od,y\in P_\ev, v\in V_\od,w\in V_\ev$. We have
    \begin{align}
        s_{\rtimes}\bigl((x+v)\cdot (y+w)\bigl)&=s\bigl(xy+\pi(x)(w)+\pi(y)(v)\bigl)\label{aaaa}\\\nonumber
        &=s(xy)+\gamma(xy)\pi(x)(w)+\gamma(xy)\pi(y)(v);\\\nonumber
        &~\\\label{bbbb}
        (y+w)^2\cdot s_{\rtimes}(x+v)&=y^2\cdot\bigl(s(x)+\gamma(x)(v)\bigl)=y^2s(x)+\pi(y^2)\gamma(x)(v);\\\nonumber
        &~\\\nonumber
        (x+v)\cdot(y+w)\cdot\{x+v,y+w\}_{\rtimes}&=\bigl(xy+\pi(x)(w)+\pi(y)(v)\bigl)\cdot\bigl(\{x,y\}+\gamma(x)(w)+\gamma(y)(v) \bigl)\nonumber\\\label{cccc}
        &=xy\{x,y\}+\pi(x)\pi(y)\gamma(x)(w)+\pi(x)\pi(y)\gamma(y)(v)\\\nonumber
        &~~\quad+\pi(\{x,y\})\pi(x)(w)+\pi(\{x,y\})\pi(y)(v).
    \end{align}
    Thus, by summing \eqref{aaaa}+\eqref{bbbb}+\eqref{cccc}, we get
    \begin{align*}
        &\gamma(xy)\pi(x)+\pi(x)\pi(y)\gamma(x)+\pi(\{x,y\})\pi(x)\\
        =~&\cancel{\pi(x)\gamma(y)\pi(x)}+\pi(y)\gamma(x)\pi(x)+\pi(x)\pi(y)\gamma(y)+\cancel{\pi(x)\gamma(y)\pi(x)}+\underset{=~0}{\underbrace{\gamma(y)\pi(x)^2}}=0,\\
    \end{align*} as maps $V_\ev\rightarrow V_\ev$, since Eq. \eqref{eq:rep3} holds. Moreover, we have
    \begin{align*}
        &\gamma(xy)\pi(y)+\pi(y^2)\gamma(x)+\pi(x)\pi(y)\gamma(y)+\pi(\{x,y\})\pi(y)\\
        =~&\cancel{\pi(x)\gamma(y)\pi(x)}+\pi(y)\gamma(x)\pi(y)+\pi(y)^2\gamma(x)+\pi(x)\pi(y)\gamma(y)\\
        &+\cancel{\pi(x)\gamma(y)\pi(x)}+\rho(y)\pi(x)\pi(y)=0,\\
    \end{align*} as maps $V_\od\rightarrow V_\ev$, since Eq. \eqref{eq:rep4} holds. Therefore, we have
   $$ s_{\rtimes}\bigl((x+v)\cdot (y+w)\bigl)= (y+w)^2\cdot s_{\rtimes}(x+v)+(x+v)\cdot(y+w)\cdot\{x+v,y+w\}_{\rtimes},$$ and the result follows.
\end{proof}

\subsection{Cohomology of Poisson superalgebras in characteristic 2}\label{sec:cohopoiss}In this Section, we introduce a cohomology theory for Poisson superalgebras in characteristic 2 that is derived from the cohomology of Lie superalgebras, see Section \ref{sec:cohoLie}.

Let $A$ be an associative supercommutative superalgebra. For $k\geq 0$, we denote by $\mathfrak{X}^k(A)$ the space of skew-symmetric $k$-derivations of $A$, that is, multilinear maps $\varphi\in\Hom(\wedge^k A,A)$ satisfying
\begin{equation}\label{eq:derivation}
   \varphi(xy,z_2,\cdots,z_k)=x\varphi(y,z_2,\cdots,z_k)+y \varphi(x,z_2,\cdots,z_k),~\forall x,y,z_2,\cdots,z_k\in A.
\end{equation}

\subsubsection{Poisson cochains} Let $(P,\{-,-\},s)$ be a Poisson superalgebra in characteristic $2$. For $n=0,1$, we set $XC^0_{\rm po}(P)=A$ and $XC^1_{\rm po}(P)=\mathfrak{X}^1(P)$. For $n\geq 2$, a pair  $(\varphi,\omega)\in XC^n(P;P)$ is a Poisson superalgebra $n$-cochain if $\varphi\in \mathfrak{X}^n(P)$ and moreover if for all $\lambda\in \K$ and all $x\in P_\od,y\in P_\ev$ and all $z_2,\cdots,z_{n-1}\in P$, we have (where  $2\leq i\leq n-1$):
\begin{align}
\label{eq:res PA1}		
\omega( x y, z_2,\cdots,z_{n-1})=&~y^2\omega(x,z_2,\cdots,z_{n-1})
+xy\varphi(x,y,z_2,\cdots,z_{n-1}),\\
\label{eq:res PA2}			
\omega( x, z_2,\cdots,z_i z'_i,\cdots,z_{n-1})=&~z_i \omega(x,z_2,\cdots,z'_i,\cdots,z_{n-1})\\\nonumber&~+z'_i\omega(x,z_2,\cdots,z_i,\cdots,z_{n-1}).
\end{align}
The space of all pairs $(\varphi,\w)\in XC^n(P;P)$ satisfying \eqref{eq:res PA1} and \eqref{eq:res PA2} is denoted by $XC^n_{\rm po}(P)$. The coboundary operators
${\rm d}_{\rm po}^n:XC^n_{\rm po}(P)\rightarrow XC^{n+1}_{\rm po}(P)$ for
$n\geq 0$ are induced by the coboundary operators for Lie superalgebras, see Section \ref{sec:diff.op.p=2} and Proposition \ref{prop:poisson-diff} below.

We denote by
$XZ_{\rm po}^n(P)=\Ker(d^n_{\rm po})$ the Poisson $n$-cocycles and $XB_{\rm po}^n(P)=\text{Im}(d^{n-1}_{\rm po})$ the Poisson $n$-coboundaries. 

Following \cite[Section 2.2]{BBE}, we can define a graduation on the cochain spaces $XC^n_{\rm po}(P)$ in the following way. Let $(\varphi,\w)\in XC^n_{\rm po}(P).$ In the case where $\varphi$ is an even map, then $\text{Im}(\w)\subset P_\ev$; whereas if $\varphi$ is odd, then $\text{Im}(\w)\subset P_\od$. Therefore, we can define a graduation on the space $XC^n_{\rm po}(P)$ by setting
$$|(\varphi,\w)|:=|\varphi|,~\forall (\varphi,\w)\in XC^n_{\rm po}(P).$$ Therefore, we have a decomposition $XC^n_{\rm po}(P)=XC^n_{\rm po}(P)_\ev\oplus XC^n_{\rm po}(P)_\od$, where
\begin{equation}\label{eq:graduation}
\begin{array}{rcll}
    XC^n_{\rm po}(P)_\ev&=&\{(\varphi,\w)\in XC^n_{\rm po}(P),~|\varphi|=\ev\};\\
     XC^n_{\rm po}(P)_\od&=&\{(\varphi,\w)\in XC^n_{\rm po}(P),~|\varphi|=\od\}.    
\end{array}
\end{equation}

\sssbegin{Proposition}\label{prop:poisson-diff}
Let $(P,\{-,-\},s)$ be a  Poisson superalgebra in characteristic $2$, and let $n\geq 0$. For all $(\varphi,\omega)\in  XC^{n}_{\rm po}(P)$, we have  $d^{n}_{\rm po}(\varphi,\omega)\in  XC^{n+1}_{\rm po}(P)$. 
\end{Proposition}
\begin{proof}
Let $x\in P_\od,~y\in P_\ev,~Z:=z_2,\cdots,z_n\in P$. We also denote by $\widehat{Z_i}:=(z_2,\cdots,\hat{z_i},\cdots,z_n)$, $\widehat{Z_{i,j}}:=(z_2,\cdots,\hat{z_i},\cdots,\hat{z_j},\cdots,z_n)$, and so on, where the hat means that the term is omitted. We have
\begin{align*}
  \delta^n\omega(xy,Z)
  =~&\{xy,\varphi(xy,Z)\}+ \displaystyle\sum_{i=2}^{n}\{z_i,\omega(xy,\widehat{Z_i})\} +\varphi(s(xy),Z)\\
  &+\displaystyle \sum_{i=2}^{n}\varphi\left(\{xy,z_i\},xy,\widehat{Z_i}\right) +\displaystyle \sum_{2\leq i<j\leq n}\omega\left(xy,\{z_i,z_j\},\widehat{Z_{i,j}}  \right)\\
=~&y^2\delta^n\omega(x,Z)+xy\,d_{\ce}^n \varphi(x,y,Z)\\
&+\{z_i,y\}y\, \omega(x,\widehat{Z_i}) + \{z_i,x\}y\, \varphi(x,y,\widehat{Z_i})+\{z_i,y\}x\, \varphi(x,y,\widehat{Z_i}    \Big)\\
=~&y^2\delta^n\omega(x,z_2,\cdots,z_n)+xy\, d_{\ce }^n \varphi(x,y,z_2,\cdots,z_n).
\end{align*}
    Moreover, we have
\[
\begin{array}{l}
\delta^n\omega(x,z_2 z'_2,\widehat{Z_2})\\[2mm]
=\{x,\varphi(x,z_2 z'_2,\widehat{Z_2})\}+\{z_2z_2',\omega(x,\widehat{Z_2})\}\\[2mm] \quad + \displaystyle \sum_{i=3}^{n}\{z_i,\omega(x,z_2z_2',\widehat{Z_{2,i}})\}+\varphi(s(x),z_2z_2',\widehat{Z_2})+\varphi\left(\{x,z_2z_2'\},x,\widehat{Z_2} \right)\\[2mm]
\quad +\displaystyle \sum_{i=3}^{n}\varphi\left(\{x,z_i\},x,z_2z_2',\widehat{Z_{2,i}} \right)
+\displaystyle 
 \sum_{i=3}^n\omega\left(x,\{z_2z_2',z_i\},\widehat{Z_{2,i}}  \right)\\[2mm]
\quad + \displaystyle \sum_{3\leq i<j\leq n}\omega\left(x,\{z_i,z_j\},z_2z_2',\widehat{Z_{2,i,j}}  \right)\\[3mm]
=z_2\delta^n\omega(x, z'_2,\widehat{Z_2})+z'_2\delta^n\omega(x, z_2,\widehat{Z_2}).
    \end{array}
    \] The result follows.
\end{proof}
Since $d_{\rm po}^n$ is induced by the coboundary operator for the Lie superalgebra $(P,\{-,-\},s)$, we have $d_{\rm po}^n(\varphi,\omega) \in C_{\rm po}^{n+1}(P)$ and $d_{\rm po}^{\,n+1} \circ d_{\rm po}^n = 0$. It follows that we obtain a well-defined cochain complex $\bigl(\bigoplus_{n \ge 0} XC_{\rm po}^n(P), d_{\rm po}^n\bigr)$, and the associated $n$-th cohomology group $XH_{\rm po}^n(P)$ is also well defined, and admits a $\Z_2$-graduation induced by \eqref{eq:graduation}.

\sssbegin{Remark} From the definition of Poisson derivations in Section \ref{def:poisson}, it follows that a linear map $\delta:P\mapsto P$ is a Poisson derivation if and only of it is a Poisson 1-cocycle.
\end{Remark}
\subsection{Formal deformations}

We briefly consider formal deformations of Poisson superalgebras as an application of the cohomology introduced in Section \ref{sec:cohopoiss}. A formal deformation of a Poisson superalgebra $(P,\{-,-\},s)$ is a Poisson superalgebra $(P[[t]],\mu_t,\w_t)$, where $P[[t]]$ is the $\K[[t]]$-algebra of formal power series with coefficients in $P$, the maps $\mu_t$ and $\w_t$ are defined by 
\begin{equation*}
    \mu_t(x,y)=\{x,y\}+\sum_{i\geq1}t^i\mu_i(x,y);\quad \w_t(z)=s(z)+\sum_{i\geq1}t^i\w_i(z),\quad\forall x,y\in P,~\forall z\in P_\od,
\end{equation*} with $(\mu_i,\w_i)\in XC^2_{\rm po}(P)_\ev.$ Note that since $(\mu_i,\w_i)\in XC^2_{\rm po}(P),$ the conditions $(i)$ and $(ii)$ of Definition \ref{def:poisson} are automatically satisfied by the maps $\mu_t$ and $\w_t$.

Two formal deformations $(P[[t]],\mu_t,\w_t)$ and $(P[[t]],\mu_t',\w_t')$ of a Poisson superalgebra\\ $(P,\{-,-\},s)$ are called \textit{equivalent} if there exists an (even) isomorphism of Poisson superalgebras 
\begin{equation}\label{eq:equiv}\psi_t:(P[[t]],\mu_t,\w_t)\to(P[[t]],\mu_t',\w_t'),\quad\psi_t=\id+\sum_{i\geq 1}t^i\psi_i,
\end{equation}
where $\psi_i:P\to P$ are linear maps.
Let $k\geq 1.$ A formal deformation of order $k$ is given by a Poisson superalgebra $(P_k[[t]],\mu_{(k)},\w_{(k)})$, where $P_k[[t]]=P[[t]]/\langle t^{k+1} \rangle$ and
\begin{equation}\label{eq:order-k}
    \mu_{(k)}(x,y)=\{x,y\}+\sum_{i\geq1}^kt^i\mu_i(x,y);\quad \w_{(k)}(z)=s(z)+\sum_{i\geq1}^kt^i\w_i(z),\quad\forall x,y\in P,~\forall z\in P_\od.
\end{equation} The formal deformations of order $k=1$ are called \textit{infinitesimal deformations}.

\sssbegin{Lemma}
Let $k\geq 1$ and let $\bigl(P_k[[t]],\mu_{(k)},\w_{(k)}\bigl)$ be a formal deformation of order $k$ of a Poisson superalgebra $\bigl(P,\cdot,\{-,-\},s\bigl)$. Then, we have $(\mu_1,\w_1)\in XZ^2_{\rm po}(P)$.
\end{Lemma}
\begin{proof}
    By expanding Eq. \eqref{eq:order-k}, we have $(\mu_1,\w_1)\in XC^2_{\rm po}(P)$. The 2-cocycle condition is covered by \cite[Theorem 6.3]{BM}.
\end{proof}

Next, we discuss the problem of extending a deformation of order $k$ to a deformation of order $k+1$. Let $(P,\{-,-\},s)$ be a Poisson superalgebra and let $(P_k[[t]],\mu_{(k)},\w_{(k)})$ be deformation of order $k$. For $k\geq 1$, we define the maps
\small{
\begin{align*}
    \obs^{(1)}_{k+1}(x,y,z)&:=\sum_{i=1}^k\Bigl(\mu_i(x,\mu_{k+1-i}(y,z))+\mu_i(y,\mu_{k+1-i}(z,x))+\mu_i(z,\mu_{k+1-i}(x,y))   \Bigl), \quad \forall~x,y,z\in P;\\
    \obs^{(2)}_{k+1}(x,y)&:=\sum_{i=1}^k\Bigl(\mu_i\bigl(y,\w_{k+1-i}(x)\bigl)+\mu_i\bigl(x,\mu_{k+1-i}(x,y)\bigl)   \Bigl),\quad \forall x\in P_\od,~\forall y\in P_\od.
\end{align*}}

\normalsize{}
\sssbegin{Lemma}\label{lem:above}
    We have $\bigl(\obs^{(1)}_{k+1},\obs^{(2)}_{k+1}\bigl)\in XC^3_{\rm po}(P)_\ev$.
\end{Lemma}
\begin{proof}
    By \cite[Theorem 6.3]{BM}, we have $\bigl(\obs^{(1)}_{k+1},\obs^{(2)}_{k+1}\bigl)\in XC^3_{\rm po}(P;P)$. Moreover, we have
    \begin{equation}
        \obs^{(1)}_{k+1}(xy,z,v)=x~\obs^{(1)}_{k+1}(y,z,v)+y~\obs^{(1)}_{k+1}(x,z,v),~\forall x,y,z,v\in P.
    \end{equation} Let $x\in P_\od$ and $y,z\in P$. We show that the pair $\bigl(\obs^{(1)}_{k+1},\obs^{(2)}_{k+1}\bigl)$ satisfies Eq. \eqref{eq:res PA1}.
\footnotesize{\begin{align*}
    &\obs^{(2)}_{k+1}(xy,z)\\=&~\sum_{i=1}^k\Bigl(\mu_i\bigl(z,\w_{k+1-i}(xy)\bigl)+\mu_i\bigl(xy,\mu_{k+1-i}(xy,z)\bigl) \\
    =&~\sum_{i=1}^k\Bigl(y^2\mu_i(z,\w_{k+1-i}(x)) +xy\cdot\mu_i(z,\mu_{k+1-i}(x,y))+\mu_{k+1-i}(x,y)\mu_i(z,xy) \Bigl)\\
    &~+\sum_{i=1}^k\Bigl(x\mu_{k+1-i}(y,z)\mu_i(y,x)+xy\cdot\mu_i(y,\mu_{k+1-i}(x,z))     \Bigl)\\
    &~+\sum_{i=1}^k\Bigl( y^2\mu_i(x,\mu_{k+1-i}(x,z))+y\mu_{k+1-i}(x,z)\mu_i(y,x)+xy\cdot\mu_i(x,\mu_{k+1-i}(y,z))     \Bigl)\\
    =&~y^2\obs^{(2)}_{k+1}(x,z)+xy\cdot\obs^{(1)}_{k+1}(x,y,z)\\
    &~+\sum_{i=1}^k\Bigl(x\mu_{k+1-i}(x,y)\mu_i(z,y)+y\mu_{k+1-i}(x,y)\mu_i(z,x)+x\mu_{k+1-i}(y,z)\mu_i(y,x)+y\mu_{k+1-i}(x,z)\mu_i(y,x)  \Bigl)\\
    =&~y^2\obs^{(2)}_{k+1}(x,z)+xy\cdot\obs^{(1)}_{k+1}(x,y,z).
\end{align*}}\normalsize{}  Moreover, Eq. \eqref{eq:res PA2} is readily checked. Therefore, we have$\bigl(\obs^{(1)}_{k+1},\obs^{(2)}_{k+1}\bigl)\in XC^3_{\rm po}(P)$.
\end{proof}

\sssbegin{Proposition}\label{prop:obs}
   Let $\bigl(P_k[[t]],\mu_{(k)},\w_{(k)}\bigl)$ be a formal deformation of order $k$ of a Poisson superalgebra $\bigl(P,\cdot,\{-,-\},s\bigl)$. Let $(\mu_{k+1},\w_{k+1})\in XC^2_{\rm po}(P)$. Then, $\bigl(P_{k+1}[[t]],\mu_{(k+1)},\w_{(k+1)}\bigl)$ is a formal deformation of order $k+1$ of $\bigl(P,\{-,-\},s\bigl)$ if and only if
     $$	\bigl( \obs^{(1)}_{k+1}, \obs^{(2)}_{k+1}\bigl)=d_{\rm po}^2\left(\mu_{k+1},~\w_{k+1} \right),$$
\end{Proposition}

     \begin{proof}
Suppose that $\bigl(P_{k+1}[[t]],\mu_{(k+1)},\w_{(k+1)}\bigl)$ is a formal deformation of order $k+1$ of $\bigl(P,\{-,-\},s\bigl)$.  It follows from \cite[Theorem 6.3]{BM} and Lemma \eqref{lem:above} that $\bigl( \obs^{(1)}_{k+1}, \obs^{(2)}_{k+1}\bigl)=d_{\rm po}^2\left(\mu_{k+1},~\w_{k+1} \right).$ The equivalence follows from the fact that we are working over $P_{k+1}[[t]]$. 
     \end{proof}

In the following Proposition, we investigate equivalent deformations and show that equivalent infinitesimal deformations are classified by the even part of the second cohomology space.
\sssbegin{Proposition}\label{prop:equiv}
    Let $k\geq 1$, and let $\bigl(P_k[[t]],\mu_{(k)},\w_{(k)}\bigl)$ and $\bigl(P_k[[t]],\mu_{(k)}',\w_{(k)'}\bigl)$ be two equivalent deformations of order $k$ of $\bigl(P,\{-,-\},s\bigl)$. Then,
    \begin{enumerate}
    \item[($i$)] we have $(\mu_1,\w_1)\equiv (\mu_1',\w_1')\in \mathrm{H}^2_{\rm po}(P)_\ev;$
    \item[($ii$)] the second cohomology space $\mathrm{H}^2_{\rm po}(P)_\ev$ classifies up to equivalence the infinitesimal deformations of $\bigl(P,\{-,-\},s\bigl)$.  
    \end{enumerate}
    \end{Proposition}
    
\begin{proof}
    Let us prove the first item. Let $\psi_t=\id+\sum_{i\geq1}t^i\psi_i$ be the formal isomorphism realizing the equivalence, see Eq. \eqref{eq:equiv}. For all $x,y\in P$, we have $$\psi_t(xy)=\psi_t(x)\psi_t(y)=xy+t\bigl(\psi_1(x)y+x\psi_1(y)\bigl)\quad \mod t^2.$$ Since $\psi_t(xy)=xy+t\psi_1(xy)\,\mod t^2$, we obtain that $\psi_1(xy)=\psi_1(x)y+x\psi_1(y)$. Therefore,  $\psi_1\in C^1_{\rm po}(P)_\ev$. Then, by \cite[Theorem 6.4]{BM}, we have $(\mu_1,\w_1)=(\mu_1',\w_1')+d^1_{\rm po}\psi_1.$ The second item is a straightforward corollary of the first one.
\end{proof}

\sssbegin{Example} Let us consider the $(1|1)$-dimensional Poisson superalgebra $P$ spanned by $e|f$ (even$|$odd) and equipped with the associative supercommutative product $ee=e,~ef=f$, the bracket $\{-,-\}=0$ and the squaring $s(f)=e$. Then, we have
$$ XH^1_{\rm po}(P;P)=\Span\{e\otimes f^*;~f\otimes f^*\};\quad XH^2_{\rm po}(P;P)=\Span\{(0,\w)\},$$ where $\w(f)=e.$
\end{Example}

\sssbegin{Example} Let us consider the $(2|2)$ dimensional Poisson superalgebra $P$ of Example \ref{ex:example}. Then, we have 
$$XH^1_{\rm po}(P;P)=\Span\{e_2\otimes e_3^*;~e_2\otimes e_4^*;~e_3\otimes e_3^*;~e_3\otimes e_4^*;~e_4\otimes e_3^*;~e_4\otimes e_4^*\};$$
$$XH^2_{\rm po}(P;P)=\Span\{(\varphi_i,0);~(0,\w_j),~i=1,2,3;~j=1,2,3,4\},$$
where
$$\varphi_1=e_2\otimes e_2^*\wedge e_4^*;~\varphi_2=e_2\otimes e_3^*\wedge e_4^*;~\varphi_3=e_3\otimes e_3^*\wedge e_4^*,$$

$$\w_1(e_3)=e_1;~\w_2(e_3)=e_2;~\w_3(e_4)=e_1;~\w_1(e_4)=e_2.$$
\end{Example}

\subsection{The universal enveloping algebra of a Poisson superalgebra}
This section aims at defining the universal enveloping algebra of a Poisson superalgebra in characteristic $p=2$. We follow the constructions of \cite{La} and show that the universal enveloping algebra can be realized as the universal enveloping algebra of an associated Lie-Rinehart superalgebra.

\subsubsection{Universal enveloping algebra}\label{def:UEApoisson} Let $(P,\{-,-\},s)$ be a Poisson superalgebra. The \textit{universal enveloping algebra} of $(P,\{-,-\},s)$ is a triple $(U(P),m,h)$, where $U$ is an associative superalgebra, $m:P\rightarrow U(P)$ is an even superalgebras morphism, $h:P\rightarrow U(P)$ is an even Lie superalgebras morphism (where $U(P)$ is seen as a Lie superalgebra with the commutator and the square power) such that
 \begin{align}
    h(xy)&=m(x)h(y)+m(y)h(x),&\forall x,y\in P;\label{eq:UEA1}\\
    m\bigl(\{x,y\}\bigl)&=m(x)h(y)+h(y)m(x),& \forall x,y\in P;\label{eq:UEA2}\\
    m(y)h(x)m(x)&=m(x)m(y)h(x),&\forall x\in P_\od,~\forall y\in P_\ev;\label{eq:UEA3}\\
    m(y)^2h(x)&=m(y)h(x)m(y)+m(x)m(y)h(y)\label{eq:UEA4}\\\nonumber&~~\quad+h(y)m(x)m(y),& \forall x\in P_\od,~\forall y\in P_\ev,
    \end{align}
    which is universal among triples $(B,m',h')$ satisfying the same properties.
Clearly stated, the universal property reads as follows. 
\begin{equation}\label{eq:UEA_univ}
\begin{minipage}{0.85\textwidth}
For any triple $(B,m',h')$ where $B$ is an associative superalgebra,
$m':P\to B$ is an even superalgebra morphism, and
$h':P\to B$ is an even Lie superalgebra morphism satisfying
conditions \eqref{eq:UEA1}--\eqref{eq:UEA4}, there exists an unique even
associative superalgebra morphism $f:U(P)\to B$ such that $f\circ h=h'$ and $f\circ m=m'$.
\end{minipage}
\end{equation}


Let $(P,\{-,-\},s)$ be a a Poisson superalgebra \textit{with unit} (for the associative multiplication). The rest of this section is devoted to construct $\bigl(U(P),m,h\bigl)$ explicitly by considering the universal enveloping algebra of the Lie-Rinehart superalgebra of K\"{a}hler differentials of $P$.

\subsubsection{Lie-Rinehart superalgebra of K\"{a}hler differentials in characteristic $2$}

   Let $A$ be an associative algebra \textit{with unit}. We denote by $S$ the free $A$-module spanned by symbols $\D a$ for all homogeneous $a\in A$. The degree is given by $|\D a|=|a|$. Let 
   $R$ be the ideal of $A$ generated by elements of the form
   $$ \D(\lambda a+\mu b)+\lambda\D a+\mu\D b,~\D(ab)+a\D b+b\D a,\text{ and }~\D\lambda,~\forall a,b\in A,~\forall \lambda,\mu\in \K.$$
   We denote by $\Omega_A:=S/R$ the \textit{module of K\"{a}hler differentials}.

The K\"{a}hler differentials satisfy the following universal property.  For every even derivation $D:A\rightarrow A$, there exists an unique $A$-linear map $\hat{D}:\Omega_A\rightarrow A$ such that $\hat{D}\circ \D=D$.  Thus, for any $A$-module $M$, precomposition by $\D$ yields a group isomorphism (see \cite[Proposition 6.1]{La})
\begin{equation}\label{eq:univ-kahler}
    \Hom_A\bigl(\Omega_A,M\bigl)\cong \Der(A,M).
\end{equation}

 Let $\bigl(P,\{-,-\},s\bigl)$ be a Poisson superalgebra in characteristic $2$ \textit{with unit}. Then, one can consider its module of K\"{a}hler differentials $\Omega_P$. In the following, we show that there is a natural Lie-Rinehart superalgebra associated to the pair $(P,\Omega_P)$, see also \cite[Theorem 3.8]{Hu} and \cite[Example 6.2]{La} for the characteristic $0$ case; \cite[Theorem 4.3.2]{BEL} and \cite[Theorem 8.2]{BYZ} for the restricted version.

\sssbegin{Lemma}
  Let $\bigl(P,\{-,-\},s\bigl)$ be a Poisson superalgebra in characteristic $2$ and let $\Omega_P$ be its module of K\"{a}hler differentials. Then $\bigl(\Omega_P,[-,-],\widetilde{s}\bigl)$ is a Lie superalgebra with the bracket
    \begin{equation}\label{kahlerbracket}
        \bigl[x\D y,z\D w\bigl]:=xz\D\{y,w\}+x\{y,z\}\D w+z\{w,x\}\D y,~\forall x\D y,z\D w\in\Omega_P;
    \end{equation}
and the squaring defined by
    \begin{equation}\label{kahlersquaring}
        \begin{cases}             
                \widetilde{s}(x\D y):=x^2\D s(y)+x\{x,y\}\D y,&~\forall x\in P_\ev,~y\in P_\od;\\
                 \widetilde{s}(x\D y):=x\{x,y\}\D y,&~\forall x\in P_\od,~y\in P_\ev.\\
        \end{cases}
    \end{equation}
\end{Lemma}

\begin{proof}
We prove the Lemma in the case where $x$ is even and $y$ is odd. Let $z\D w\in \Omega_P$. We have
\begin{align}\label{eq:tr1}
    \bigl[\widetilde{s}(x\D y),zdw\bigl]=&~x^2zd\{s(y),w\}+x^2\{s(y),z\}\D w+xz\{x,y\}\D \{y,w\}\\\nonumber&+x\{x,y\}\{y,z\}\D w+xz\{w,\{x,y\}\D y+z\{x,y\}\{w,x\}\D y.
\end{align} Moreover, we have
\begin{align}
    \bigl[x\D y,xz\D \{y,w\}\bigl]=&~x^2z\D \{y,\{y,w\}\}+x^2\{y,z\}\D \{y,w\}\\\nonumber&+xz\{y,x\}\D \{y,w\}+xz\{\{y,w\},x\}\D y;
\end{align}
\begin{align}
    \bigl[x\D y,x\{y,z\}\D w\bigl]=&~x^2\{y,z\}\D \{y,w\}+x^2\{y,\{y,z\}\D w\\\nonumber&+x\{y,z\}\{y,x\}\D w+x\{y,z\}\{w,x\}\D y;
\end{align}
\begin{align}\label{eq:tr4}
    \bigl[x\D y,z\{w,z\}\D y\bigl]=&~xz\{y,\{w,x\}\}\D y+x\{y,z\}\{w,x\}\D y+z\{x,y\}\{w,x\}\D y.
\end{align}
Comparing Eqs. \eqref{eq:tr1}-\eqref{eq:tr4}, we have
$$\bigl[\widetilde{s}(x\D y),z\D w\bigl]=\bigl[x\D y,\bigl[x\D y,z\D w\bigl]\bigl].$$ Thus, it follows from Proposition \ref{prop:jacobson} that $\bigl(\Omega_P,[-,-],\Tilde{s}\bigl)$ is a Lie superalgebra.
\end{proof}

\sssbegin{Proposition}\label{prop:kahler-LR}
     Let $\bigl(P,\{-,-\},s\bigl)$ be a Poisson superalgebra in characteristic $2$ and let $\Omega_P$ be its module of K\"{a}hler differentials viewed as a Lie superalgebra with the bracket \eqref{kahlerbracket} and the squaring \eqref{kahlersquaring} and endowed with the Lie superalgebras map 
     \begin{equation}
          \rho:\Omega_P\rightarrow \Der_{\rm{asso}}(P),~x\D y\mapsto x\{y,-\}.
     \end{equation}
Then, $(P,\Omega_P,\rho)$ is a Lie-Rinehart superalgebra. 
\end{Proposition}

\begin{proof}
Let $x,y,z\in P$. There are four cases to consider, depending on the parity of the elements.

\underline{The case where $y$ is odd and $x,z$ are even.} In that case, we have
\begin{align*}
    \widetilde{s}(z\cdot x\D y)&=(zx)^2\D s(y)+zx\{zx,y\}\D y\\
            &=z^2x^2\D s(y)+z^2x\{x,y\}\D y+zx^2\{z,y\}\D y\\
            &=z^2\widetilde{s}(x\D y)+\rho(z\cdot x\D y)(z)x\D y.
\end{align*}

\underline{The case where $x$ is odd and $y,z$ are even.} In that case, we have
$$\widetilde{s}(z\cdot x\D y)=z^2x\{x,y\}\D y=z^2\widetilde{s}(x\D y)+\rho(z\cdot x\D y)(z)x\D y.$$

\underline{The case where $z$ is odd and $x,y$ are even.} In that case, we have
$$\widetilde{s}(z\cdot x\D y)=z^2x\{x,y\}\D y=\rho(z\cdot x\D y)(z)x\D y.$$

\underline{The case where $x,y,z$ are odd.} In that case, both $\widetilde{s}(z\cdot x\D y)$ and $\rho(z\cdot x\D y)(z)x\D y$ vanish.
The conclusion follows from Lemma \ref{lem:ext-rinhart}.
\end{proof}
The following Theorem shows that the universal enveloping algebra of a Poisson superalgebra $P$ can be realized as the Lie-Rinehart enveloping algebra of the Lie-Rinehart superalgebra $(P,\Omega_P,\rho)$ described in Proposition \ref{prop:kahler-LR}, (see \cite[Proposition 6.3]{La} for the characteristic zero case).
\sssbegin{Theorem}\label{thm:univ-iso}
 Let $\bigl(P,\{-,-\},s\bigl)$ be a Poisson superalgebra in characteristic $2$ and let $\Omega_P$ be its module of K\"{a}hler differentials. Then, the universal enveloping algebra $U(P,\Omega_P)$ of the Lie-Rinehart superalgebra $(P,\Omega_P,\rho)$ satisfies the universal property \eqref{eq:UEA_univ}.
\end{Theorem}

\begin{proof}
    Denote by $i_P:P\rightarrow U(P,\Omega_P)$ and $i_\Omega:\Omega_P\rightarrow U(P,\Omega_P)$ the canonical inclusions and let $j=i_\Omega\circ \D$. We will show that the triple $\bigl(U(P,\Omega_P),i_\Omega,j\bigl)$ satisfy the universal property \eqref{eq:UEA_univ}. It follows from the construction of $U(P,\Omega_P)$ that $i_P$ and $j$ are even and that $i_p$ is an associative superalgebras morphism. Let us show that $j$ is a morphism of Lie superalgebras. Let $x\in P_\od$. we have
    $$j\bigl(s(x)\bigl)=i_{\Omega}\circ\widetilde{s}(\D x)=i_{\Omega}(\D x)^2=j(x)^2.$$ We also have $j\bigl(\{x,y\}\bigl)=j(x)j(y)+j(y)j(x)$ for any $x,y\in P$, see \cite[Proposition 6.3]{La}. Thus, the map $j$ is a morphism of Lie superalgebras. Similarly to \cite[Proposition 6.3]{La}, we have
    \begin{equation}
    j(xy)=i_P(x)j(y)+i_P(y)j(x),\quad\forall x,y \in P,
    \end{equation}
    and
     \begin{equation}\label{eq:qwe}
     i_P\bigl(\{x,y\}\bigl)=j(x)i_P(y)+i_P(y)j(x),\quad\forall x,y \in P.
     \end{equation}
    Note that taking $x=y$ in Eq.\eqref{eq:qwe} yields $j(x)i_P(x)=i_P(x)j(x),~\forall x\in P$.
    Let $x\in P_\od$ and $y\in P_\ev$. Let us prove that $i_P$ and $j$ satisfy Eq. \eqref{eq:UEA3}. We have
    $$i_P(y)j(x)i_P(x)=i_P(y)i_P(x)j(x)=i_P(x)i_P(y)j(x),$$ therefore Eq. \eqref{eq:UEA3} is satisfied. Finally, let us consider Eq. \eqref{eq:UEA4}. Notice that $\forall x,y,z\in P$, we have
    \begin{align*}
        \rho(\D x)(y)\D z&=[\D x,y\D z]+y[\D x,\D z]=2y\D\{x,z\}+\{x,y\}\D z=\{x,y\}\D z.                     
    \end{align*} It follows that $\rho(\D x)(y)\D z=\{x,y\}\D z=\rho(\D y)(x)\D z.,$ $\forall x,y,z\in P$.
    Let $x\in P_\od$ and $y\in P_\ev$. We have
    \begin{equation*}
        i_P(y)i_P(y)j(x)=i_P(y)\bigl(j(x)i_p(y)+i_P(\rho(\D x)(y))\bigl)=i_P(y)j(x)i_p(y)+i_P(y)i_P\bigl(\rho(\D x)(y)\bigl)
    \end{equation*} and
    \begin{align*}
        i_P(x)i_P(y)j(y)&=i_P(y)i_P(x)j(y)\\
        &=i_P(y)j(y)i_P(x)+i_P(y)i_P\bigl(\rho(\D y)(x)\bigl)\\
        &=j(y)i_P(y)i_P(x)+i_P(y)i_P\bigl(\rho(\D y)(x)\bigl).
    \end{align*}
    It follows that
    \begin{align*}
        i_P(y)i_P(y)j(x)+ i_P(x)i_P(y)j(y)=~&i_P(y)j(x)i_p(y)+i_P(y)i_P\bigl(\rho(\D x)(y)\bigl)\\
        &+j(y)i_P(y)i_P(x)+i_P(y)i_P\bigl(\rho(\D y)(x)\bigl)\\
        =i_P(y)j(x)i_p(y)+j(y)i_P(y)i_P(x),
    \end{align*} thus Eq. \eqref{eq:UEA4} is satisfied. It remains to show that the triple $\bigl(U(P,\Omega_P),i_\Omega,j\bigl)$ is universal among all triples $(B,m,h)$ satisfying conditions \eqref{eq:UEA1} to \eqref{eq:UEA4}. Following \cite[Proposition 6.3]{La}, $B$ can be seen as a $P$-module with the action $x\cdot b=m(x)b,~\forall x\in P,~\forall b\in B$. Therefore, we have that $h:P\rightarrow B$ satisfies
    $$h(xy)=m(x)h(y)+m(y)h(x)=x\cdot h(y)+y\cdot h(x),$$ thus $h$ is a derivation with values in $B$. By the universal property of the K\"{a}hler differentials \eqref{eq:univ-kahler}, there exists an unique (even) $P$-linear map $\widehat{h}:\Omega_P\rightarrow B$ such that $h=\widehat{h}\circ\D$. Then, the triple $(B,m,\widehat{h})$ satisfy the following properties:
    \begin{enumerate}{}
    \item[$\textup{(}i\textup{)}$] $m:P\rightarrow B$ is an even superalgebras morphism,
      \item[$\textup{(}ii\textup{)}$] $\widehat{h}
      :\Omega_P\rightarrow B$ is an even Lie superalgebras morphism,
      \item[$\textup{(}iii\textup{)}$] we have $m\circ\rho(\D x)(y)=\widehat{h}(\D x)m(y)+m(y)\widehat{h}(\D x),~\forall x,y\in P$, and
      \item[$\textup{(}iv\textup{)}$] we have $\widehat{h}(y\D x)=m(y)\widehat{h}(\D x)~\forall x,y\in P$.
    \end{enumerate} 
    
    Therefore, by the universal property of the enveloping algebra of a Lie-Rinehart superalgebras, there is an unique map $\theta:U(P,\Omega_P)\rightarrow B$ such that $\theta\circ i_P=m$ and $\theta\circ j=\widehat{h}\D=h$. It follows that the triple $\bigl(U(P,\Omega_P),i_\Omega,j\bigl)$ satisfies the universal property \eqref{eq:UEA_univ}.
\end{proof}
\appendix
\renewcommand{\thesection}{\Alph{section}}
\renewcommand{\thesubsection}{\thesection.\arabic{subsection}}
\section{Appendix: pre-Poisson superalgebras}\label{sec:5}
In this appendix, we initiate the study ofpre-Poisson superalgebras in characteristic 2, following \cite{A}.
\subsection{Zinbiel superalgebras}
Following \cite{Lo,CFKN}, a superspace $Z=Z_\ev\oplus Z_\od$ is called (right) Zinbiel superalgebra if it is equipped with a bilinear product $*$ satisfying
\begin{align}\label{eq:Zinbiel}
    x*(y*z)&=(y*x)*z+(x*y)*z,\quad \forall x,y,z\in Z.
\end{align}
Note that by taking $x=y$ in Eq. \eqref{eq:Zinbiel}, we have $x*(x*z)=2(x*x)*z=0,\quad \forall x,z\in Z.$

\sssbegin{Proposition}\label{prop:zinbielasso}
Let $(Z,*)$ be a right Zinbiel superalgebra. Then
\begin{enumerate}
    \item[($i$)] $(Z,*)$ is left-commutative, that is, $x*(y*z)=y*(x*z),\quad \forall x,y,z\in Z$;
    \item[($ii$)] the product $x\cdot y:=x*y+y*x,~\forall x,y\in Z,$ is associative supercommutative.
\end{enumerate}
\end{Proposition}

\begin{proof} Left-commutativity is a straightforward consequence of Eq. \eqref{eq:Zinbiel}.
 Following Loday's proof (see \cite{Lo}) for all $x,y,z\in Z$, we have
\begin{align*}
    x\cdot(y\cdot z)&=(y*x)*z+(x*y)*z\\
    &+(x*z)*y+(x*z)*y\\
    &+(y*z)*x+(z*y)*x=(x\cdot y)\cdot z.
\end{align*} In particular, if $x=y$, we obtain $x\cdot(x\cdot z)=0=(x\cdot x)\cdot z,\quad\forall x,z\in Z.$
\end{proof}

\sssbegin{Theorem}[\text{\cite[Theorem 5.11]{CFKN}}]
Any finite-dimensional Zinbiel superalgebra $(Z,*)$ over any field is nilpotent.
\end{Theorem}

\subsection{Left-symmetric superalgebras in characteristic 2}\label{prelie} Following \cite{BBE}, a vector superspace $V=V_\ev\oplus V_\od$ is called a left-symmetric superalgebra $(V,\trr)$ in characteristic $p=2$ if it is endowed with a bilinear product $\trr:V\times V\rightarrow V$ satisfying
\begin{equation*}
\begin{array}{lrlll}
(i) &x\trr (y\trr z)+(x\trr y)\trr z& = & y\trr (x\trr z)+(y\trr x)\trr z,&\forall x,y,z\in V; \\[2mm]
(ii) &x\trr(x\trr y)&=&(x\trr x)\trr y,&\forall x\in V_\od,~\forall y\in V.\\[2mm]
\end{array}
\end{equation*}


\sssbegin{Proposition}[\text{\cite[Proposition 2.4.1]{BBE}}]\label{prop:LSSA-Lie}
    Let $(V,\trr)$ be a left-symmetric superalgebra. Then, $\bigl(\fg(V),[-,-],s\bigl)$ is a Lie superalgebra with $\fg(V)=V$ as superspaces and
   \begin{equation}\label{LSSA-Lie1} 
   \begin{array}{rcll}
        [x,y]&:=&x\trr y+y\trr x,~&\forall x\in V_\ev,\forall y\in V;\\
        s(x)&:=&x\trr x,~&\forall x\in V_\od.
    \end{array}
    \end{equation}
    A left-symmetric product $\trr$ on a Lie superalgebra $(V,[-,-],s)$ is called compatible with the Lie superalgebra structure if conditions \eqref{LSSA-Lie1} are satisfied.
\end{Proposition}

\subsection{Pre-Poisson superalgebras in characteristic 2} Following \cite{A}, a triple\\ $(P=P_\ev\oplus P_\od,*,\trr)$ is called pre-Poisson superalgebra if
\begin{enumerate}
    \item[($i$)] $(P,*)$ is a right Zinbiel superalgebra;
    \item[($ii$)] $(P,\trr)$ is a left-symmetric superalgebra;
    \item[($iii$)] the following compatibility conditions are satisfied for all $x,y,z\in P$:
    \begin{align}
        (x\trr y+y\trr x)*z&=x\trr (y*z)+y*(x\trr z),\quad \forall x,y,z\in P;\label{eq:preP1}\\
        (x*y+y*x)\trr z&= x*(y\trr z)+y*(x\trr z)\quad \forall x,y,z\in P.\label{eq:preP2}
    \end{align}
\end{enumerate}
Note that Eq. \eqref{eq:preP1} implies that $x\trr(x*z)=x*(x\trr z),\quad \forall x,z\in P$.

\sssbegin{Proposition}
    Let $(P,*,\trr)$ be a pre-Poisson superalgebra. Then, $(P,\cdot,[-,-],s)$ is a Poisson superalgebra, with the associative supercommutative structure $(\cdot)$ defined in Prop. \ref{prop:zinbielasso} and the Lie structure $([-,-],s)$ defined in Prop. \ref{prop:LSSA-Lie}.
\end{Proposition}

\begin{proof}
    It is well known that the bracket $[-,-]$ induces a derivation of the product $\cdot$, see e.g. \cite{A}. It remains to check Condition $(ii)$ of Definition \ref{def:poisson}. Let $x\in P_\od,~y\in P_\ev$. We have
    \begin{align*}
        s(x\cdot y)=&~(x*y)\trr(x*y)+(x*y)\trr(y*x)+(y*x)\trr(x*y)+(y*x)\trr(y*x)\\
        =&~x*(y\trr(x*y))+y*(x\trr(x*y))+x*(y\trr(y*x))+y*(x\trr(y*x)),
    \end{align*}
    Furthermore, we have
    \begin{align*}
        (x\cdot y)\cdot[x,y]=&~(x*y)*(x\trr y)+(x*y)*(y\trr x)+(y*x)*(x\trr y)+(y*x)*(y\trr x)\\
        &+(x\trr y)*(x*y)+(x\trr y)*(y*x)+(y\trr x)*(x*y)+(y\trr x)*(y*x).
    \end{align*}
    We see that some terms appear in both expressions. Since $(y\cdot y)\cdot s(x)=0,$ it remains to show that
    \begin{equation}\label{eq:step}
        x*(y\trr(x*y))=x*(y*(x\trr y))+x\trr(y*(x*y)).
    \end{equation}
    Since
    \begin{align*}
        x\trr(y*(x*y))&=x\trr(x*(y*y))\\
        &=x\trr((x*y)*y)+x\trr((y*x)*y)\\
        &=x\trr(x*(y*y))\\
        &=x*(x\trr(y*y)),
    \end{align*}
    the right-hand side of Eq. \eqref{eq:step} is equal to $$x*\bigl(y*(x\trr y)+x\trr(y*y) \bigl)=x*\bigl((x\trr y)*y+(y\trr x)*y \bigl).$$ For the left-hand side of Eq. \eqref{eq:step}, we have
    \begin{align*}
        x*(y\trr(x*y))&=x*\bigl(x*(y\trr y)+(y\trr x)*y+(x\trr y)*y\bigl)\\
        &=x*\bigl((y\trr x)*y+(x\trr y)*y\bigl).
    \end{align*} Thus, Eq. \eqref{eq:step} is satisfied and therefore, Condition $(ii)$ of Definition \ref{def:poisson} as well.
\end{proof}

Recall that an even linear operator $R$ on an associative superalgebra $(A,\cdot)$ is called \textit{Rota-Baxter operator} (see \cite{Ro}) if it satisfies
\begin{equation}
    R(x)\cdot R(y)=R\bigl(R(x)\cdot y+x\cdot R(y)\bigl),\quad \forall x,y\in A,
\end{equation}
while an even operator $Q$ on a Lie superalgebra $(L,[-,-],s)$ is called \textit{Rota-Baxter operator} (see \cite{BBE}) if it satisfies
 \begin{align}
        [Q(x),Q(y)]&=Q\bigl([Q(x),y]+[x, Q(y)]\bigl),~\forall x,y\in L;\\
        s\circ Q(x)&=Q\bigl([Q(x),x]),~\quad\forall x\in L_\od.
    \end{align}

\sssbegin{Proposition}
    Let $(P,\cdot,[-,-],s)$ be a Poisson superalgebra and let $R$ be a Rota-Baxter operator with respect to the product $\cdot$ and to the Lie structure $\bigl([-,-],s\bigl)$ simultaneously. Define new products $*$ and $\trr$ on $P$ by
    \begin{align*}
        x*y&=R(x)\cdot y;\quad x\trr y=\bigl[R(x),y\bigl],\quad\forall x,y\in P.
    \end{align*}
    Then, $(P,*)$ is a right Zinbiel superalgebra, $(P,\trr)$ is a left-symmetric superalgebra and $(P,*,\trr)$ is a pre-Poisson superalgebra.
\end{Proposition}

\begin{proof}
    Straightforward consequences of \cite[Proposition 2.4.11]{BBE} and \cite[Propositions 5.1 and 5.2]{A}.
\end{proof}

\noindent\textbf{Acknowledgement.} I thank S. Bouarroudj for his support, help, and many fruitful discussions.

\end{document}